\documentclass[a4paper,reqno,12pt]{amsart} 
\usepackage{newclude}
\usepackage[margin=1.2in]{geometry}
\usepackage{indentfirst}
\usepackage{amsfonts}
\usepackage{amscd}
\usepackage{amsthm}
\usepackage{amssymb}
\usepackage{amsmath}
\usepackage{bm}
\usepackage{thmtools}
\usepackage{mathrsfs}
\usepackage{epsfig}
\usepackage{graphicx}
\usepackage{slashed}
\usepackage[all,cmtip]{xy}
\usepackage{tikz}
\usepackage{tikz-cd}
\usepackage[textsize=footnotesize]{todonotes}
\usepackage{extpfeil}
\usetikzlibrary{decorations.markings,shapes.geometric}
\usepackage{enumitem}
\usepackage{setspace}
\usepackage{upgreek}
\usepackage{verbatim}
\usepackage[bookmarksdepth=2,linktoc=page,colorlinks,linkcolor={red!80!black},citecolor={red!80!black},urlcolor={blue!80!black}]{hyperref}

\usetikzlibrary{matrix,arrows,positioning}

\def\itemNum$#1${\item $\displaystyle#1$
   \hfill\refstepcounter{equation}(\theequation)}

\makeatletter
\providecommand\@dotsep{5}
\renewcommand{\listoftodos}[1][\@todonotes@todolistname]{%
  \@starttoc{tdo}{#1}}
\makeatother

\newtheorem{Lem}{Lemma}[section]
\newtheorem{Prop}[Lem]{Proposition}

\newtheorem*{Def}{Definition}

\theoremstyle{plain}
\newtheorem{Thm}[Lem]{Theorem}

\newtheorem{Cor}[Lem]{Corollary}

\newtheorem{thm}{Theorem}

\theoremstyle{definition}
\declaretheorem[numbered=no,name=Example,qed={\lower-0.3ex\hbox{$\triangleleft$}}]{Ex}
\newtheorem{ex}[Lem]{Example}
\newtheorem{Rem}[Lem]{Remark}
\newtheorem{Rems}[Lem]{Remarks}

\DeclareMathOperator{\Hom}{Hom}

\DeclareMathOperator{\End}{End}

\DeclareMathOperator{\Spec}{Spec}

\DeclareMathOperator{\Vect}{\mathsf{Vect}}

\DeclareMathOperator{\Modu}{\mhyphen\mathsf{Mod}}
\DeclareMathOperator{\proj}{\mhyphen\mathsf{proj}}

\newcommand{\im}{\textup{im}}
\newcommand{\coker}{\text{\textnormal{coker}}}
\mathchardef\mhyphen="2D
\newcommand{\git}{/\!\!/}

\newcommand{\op}{\text{\textnormal{op}}}
\newcommand{\id}{\text{\textnormal{id}}}

\newcommand{\Z}{\mathbb{Z}}

\newcommand{\Fun}{{\mathbb{F}_1}}

\newcommand{\ad}{\textup{ad}}

\makeatletter
\newbox\xrat@below
\newbox\xrat@above

\newcommand{\xrightarrowtail}[2][]{%
  \setbox\xrat@below=\hbox{\ensuremath{\scriptstyle #1}}%
  \setbox\xrat@above=\hbox{\ensuremath{\scriptstyle #2}}%
  \pgfmathsetlengthmacro{\xrat@len}{max(\wd\xrat@below,\wd\xrat@above)+.6em}%
  \mathrel{\tikz [>->,baseline=-.75ex]
                 \draw (0,0) -- node[below=-2pt] {\box\xrat@below}
                                node[above=-2pt] {\box\xrat@above}
                       (1.2,0) ;}}
\newcommand{\xxrightarrowtail}[2][]{%
  \setbox\xrat@below=\hbox{\ensuremath{\scriptstyle #1}}%
  \setbox\xrat@above=\hbox{\ensuremath{\scriptstyle #2}}%
  \pgfmathsetlengthmacro{\xrat@len}{max(\wd\xrat@below,\wd\xrat@above)+.6em}%
  \mathrel{\tikz [>->,baseline=-.75ex]
                 \draw (0,0) -- node[below=-2pt] {\box\xrat@below}
                                node[above=-2pt] {\box\xrat@above}
                       (2.4,0) ;}}
\renewcommand{\xtwoheadrightarrow}[2][]{%
  \setbox\xrat@below=\hbox{\ensuremath{\scriptstyle #1}}%
  \setbox\xrat@above=\hbox{\ensuremath{\scriptstyle #2}}%
  \pgfmathsetlengthmacro{\xrat@len}{max(\wd\xrat@below,\wd\xrat@above)+.6em}%
  \mathrel{\tikz [->>,baseline=-.75ex]
                 \draw (0,0) -- node[below=-2pt] {\box\xrat@below}
                                node[above=-2pt] {\box\xrat@above}
                       (1.2,0) ;}}
\newcommand{\yrightarrowtail}[2][]{%
  \setbox\xrat@below=\hbox{\ensuremath{\scriptstyle #1}}%
  \setbox\xrat@above=\hbox{\ensuremath{\scriptstyle #2}}%
  \pgfmathsetlengthmacro{\xrat@len}{max(\wd\xrat@below,\wd\xrat@above)+.6em}%
  \mathrel{\tikz [>->,baseline=-.75ex]
                 \draw (0,0) -- node[below=-2pt] {\box\xrat@below}
                                node[above=-2pt] {\box\xrat@above}
                       (0.7,0) ;}}
\newcommand{\ytwoheadrightarrow}[2][]{%
  \setbox\xrat@below=\hbox{\ensuremath{\scriptstyle #1}}%
  \setbox\xrat@above=\hbox{\ensuremath{\scriptstyle #2}}%
  \pgfmathsetlengthmacro{\xrat@len}{max(\wd\xrat@below,\wd\xrat@above)+.6em}%
  \mathrel{\tikz [->>,baseline=-.75ex]
                 \draw (0,0) -- node[below=-2pt] {\box\xrat@below}
                                node[above=-2pt] {\box\xrat@above}
                       (0.7,0) ;}}
\makeatother

\title[$K$-theory and $GW$-theory of proto-exact categories]{Group completion in the $K$-theory and Grothendieck--Witt theory of proto-exact categories}

\author[J.\,N. Eberhardt]{Jens Niklas Eberhardt}
\address{Max Planck Institute for Mathematics\\
Vivatsgasse 7\\
53111 Bonn, Germany}
\email{mail@jenseberhardt.com}

\author[O. Lorscheid]{Oliver Lorscheid}
\address{Instituto Nacional de Matem\'{a}tica Pura e Aplicada, Rio de Janeiro, Brazil}
\email{oliver@impa.br}

\author[M.\,B. Young]{Matthew B. Young}
\address{Department of Mathematics and Statistics \\ Utah State University\\
Logan, Utah 84322 \\ USA}
\email{matthew.young@usu.edu}

\date{\today}


\begin{document}

\begin{abstract}
We study the algebraic $K$-theory and Grothendieck--Witt theory of proto-exact categories, with a particular focus on classes of examples of $\mathbb{F}_1$-linear nature. Our main results are analogues of theorems of Quillen and Schlichting, relating the $K$-theory or Grothendieck--Witt theories of proto-exact categories defined using the (hermitian) $Q$-construction and group completion.
\end{abstract}

\maketitle

\begin{small} \tableofcontents \end{small}

\setcounter{footnote}{0}

\section*{Introduction}
\addtocontents{toc}{\protect\setcounter{tocdepth}{1}}

The algebraic $K$-theory of an exact category can be defined in essentially two different ways. The first way, which includes Quillen's $Q$-construction \cite{quillen1973} and Waldhausen's $\mathcal{S}_{\bullet}$-construction \cite{waldhausen1985}, uses the exact structure of the category while the second way, which includes group completion and the $+$-construction \cite{grayson1976}, uses only the underlying additive structure of the category. Each approach has its own strengths, both computational and theoretical. A celebrated theorem of Quillen \cite{grayson1976}, henceforth referred to as the Group Completion Theorem, asserts that the $K$-theories defined by the $Q$- or $\mathcal{S}_{\bullet}$-constructions and by group completion agree for split exact categories, that is, for those exact categories in which all short exact sequences are split. A closely related result, Quillen's `$Q=+$' Theorem, asserts the analogous statement for the $+$-construction under additional hypotheses on the category \cite{grayson1976}. Quillen's Group Completion Theorem was generalized to Grothendieck--Witt theory, also called hermitian $K$-theory, in Schlichting's seminal work \cite{schlichting2004}, \cite{schlichting2017}, \cite{schlichting2019}.

In this paper, we investigate analogues of Quillen's Group Completion Theorem in \emph{non-additive} contexts. To do so, in place of exact categories we study \emph{proto-exact} categories with an \emph{exact direct sum}. This is a minimal set-up in which the Group Completion Theorem can be formulated. A proto-exact category, as introduced by Dyckerhoff and Kapranov \cite{dyckerhoff2019}, is a category $\mathcal{A}$ together with two distinguished classes of morphisms, called inflations and deflations, satisfying axioms which allow for a direct adaptation of the $Q$-construction and $\mathcal{S}_{\bullet}$-construction. An \emph{exact direct sum} is a symmetric monoidal structure $\oplus$ on $\mathcal{A}$ which is compatible with the proto-exact structure but fulfills only a subset of the axioms of a categorical product or coproduct.

We focus our attention on \emph{uniquely split} proto-exact categories, that is, proto-exact categories with exact direct sum in which each short exact sequence splits in a unique way. Our main results are comparison theorems for $K$-theory and Grothendieck--Witt spaces of such categories defined using the $Q$-construction and group completion. For $K$-theory, our result is a direct analogue of Quillen's Group Completion Theorem. The situation is more subtle for Grothendieck--Witt theory and the results differ substantially from the exact case.

The main motivation for studying uniquely split proto-exact categories is their central role in $\mathbb{F}_1$-geometry \cite{deitmar2006}, \cite{connes2010}, \cite{chu2012}. A characteristic feature of $\mathbb{F}_1$-geometry is the lack of additivity in its natural constructions. For example, let $A \Modu$ be the category of modules over a commutative pointed monoid $A$ or, in geometric terms, the category of quasi-coherent sheaves over the affine monoid scheme $\Spec(A)$. The category $A \Modu$ does not have an additive structure and so, in particular, has no exact structure. However, $A\Modu$ can be given the structure of a proto-exact category with exact direct sum given by disjoint union. This notion is closely related to other partially additive structures, such as quasi-exact \cite{deitmar2006} or belian categories \cite{deitmar2012}, but has the distinct advantage of being self-dual and so well-adapted to Grothendieck--Witt theory.  

A second characteristic feature of $\mathbb{F}_1$-geometry is the distinguished role played by generalized monomial matrices, that is, matrices with at most one non-zero entry in each row and column. In our categorical context, this is reflected in the fact that exact sequences in many proto-exact categories of interest not only split, but do so in a \emph{unique} way. For example, the category $A \proj$ of finitely generated projective $A$-modules is uniquely split proto-exact, while $A \Modu$ is in general non-split.

After establishing required categorical background in Section \ref{sec:prelim}, we turn in Section \ref{sec:algKProtoEx} to $K$-theory. Let $\mathcal{A}$ be a proto-exact category with exact direct sum $\oplus$. Denote by $\mathcal{K}(A)$ and $\mathcal{K}^{\oplus}(\mathcal{A})$ the $K$-theory spaces of $\mathcal{A}$ defined via the $Q$-construction and group completion, respectively. For group completion, we regard $\mathcal{A}$ as a symmetric monoidal category by forgetting its proto-exact structure. The homotopy groups of these spaces are the corresponding algebraic $K$-theory groups $K_i(\mathcal{A})$ and $K_i^{\oplus}(\mathcal{A})$. In the setting of split exact categories, Quillen's Group Completion Theorem \cite{grayson1976} states that there is a homotopy equivalence $\mathcal{K}(\mathcal{A}) \simeq \mathcal{K}^{\oplus}(\mathcal{A})$. The first main result of this paper is a Group Completion Theorem for uniquely split proto-exact categories.

\begin{thm}[{Theorem \ref{thm:KComparison}}]
Let $\mathcal{A}$ be a uniquely split proto-exact category. Then there is a homotopy equivalence $\mathcal{K}(\mathcal{A}) \simeq \mathcal{K}^{\oplus}(\mathcal{A})$.
\end{thm}

The proof of Theorem \ref{thm:KComparison} is a modification of the original argument \cite{grayson1976}, simplified in view of the uniquely split assumption. Many of the techniques of the proof are used in the subsequent considerations for Grothendieck--Witt theory.

In Section \ref{sec:GWProtoEx} we turn to Grothendieck--Witt theory, where the situation is more delicate. Suppose that $\mathcal{A}$ has a compatible duality structure, that is, an exact functor $P: \mathcal{A}^{\op} \rightarrow \mathcal{A}$ together with coherence data exhibiting it as an involution.  The Grothendieck--Witt space $\mathcal{GW}(\mathcal{A})$ is the homotopy fibre of the forgetful functor $BQ_h(\mathcal{A}) \rightarrow BQ(\mathcal{A})$ over the zero object $0$, where $Q_h(\mathcal{A})$ is the proto-exact hermitian $Q$-construction of $\mathcal{A}$ \cite{charney1986}, \cite{uridia1990}. Again, by forgetting the proto-exact structure, one can also form a direct sum Grothendieck--Witt space $\mathcal{GW}^{\oplus}(\mathcal{A})$ using group completion. Under the assumption that $\mathcal{A}$ is split exact, the Group Completion Theorem for Grothendieck--Witt theory was proved only recently, giving a homotopy equivalence between $\mathcal{GW}(\mathcal{A})$ and $\mathcal{GW}^{\oplus}(\mathcal{A})$ \cite{schlichting2004}, \cite{schlichting2017}, \cite{schlichting2019}. The proof is considerably more involved than the corresponding result for $K$-theory. 

In the setting of Grothendieck--Witt theory of uniquely split proto-exact categories, we find that the Group Completion Theorem fails even in the basic case of $\Fun$-vector spaces (cf.\ Example \ref{ex: GW-theory of Vect F1}). However, it does hold for hyperbolic Grothendieck--Witt theory of uniquely split proto-exact categories. To state this, let $Q_H(\mathcal{A}) \subset Q_{h}(\mathcal{A})$ be the full subcategory on hyperbolic symmetric forms and write $\mathcal{GW}_H(\mathcal{A})$ for the homotopy fibre of $BQ_H(\mathcal{A}) \rightarrow BQ(\mathcal{A})$ over $0$. Let also $\mathcal{GW}_H^{\oplus}(\mathcal{A})$ be the space obtained by applying group completion to the monoidal groupoid of hyperbolic symmetric forms in $\mathcal{A}$. 

\begin{thm}[{Theorem \ref{thm:GWComparison}}]
Let $\mathcal{A}$ be a uniquely split proto-exact category with duality. Then there is a weak homotopy equivalence $\mathcal{GW}_H(\mathcal{A}) \simeq \mathcal{GW}_H^{\oplus}(\mathcal{A})$.
\end{thm}

We use Theorem \ref{thm:GWComparison} to compute the higher Grothendieck--Witt groups of $\mathcal{A}$ in terms of the symmetric monoidal groupoid of hyperbolic forms. Namely, we prove that there are abelian group isomorphisms
\[
GW_i(\mathcal{A}) 
\simeq
GW^{\oplus}_{H,i}(\mathcal{A}),
\qquad
i \geq 1,
\]
cf.\ Corollary \ref{cor:higherGWGroupsQ+}. In concrete examples, $\mathcal{GW}_H^{\oplus}(\mathcal{A})$ and $\mathcal{GW}^{\oplus}(\mathcal{A})$ can often be described using the $+$-construction, thereby allowing for computations. Such examples are discussed in the companion paper \cite{eberhardtLorscheidYoung2020b}. 

In general, the description of the homotopy type of the full space $\mathcal{GW}(\mathcal{A})$ is more complicated than that of $\mathcal{GW}_H(\mathcal{A})$ and $\mathcal{GW}^{\oplus}(\mathcal{A})$.
Under the additional assumption that $\mathcal{A}$ is combinatorial, a hypothesis that is typically satisfied in $\mathbb{F}_1$-linear contexts, we are able to describe $\mathcal{GW}(\mathcal{A})$. Here combinatorial means that any subobject $U$ of an exact direct sum $X \oplus Y$ splits into objects $U \cap X \subset X$ and $U \cap Y \subset Y$. 
\begin{thm}[{Theorem \ref{thm:GWComputation}}]
Let $\mathcal{A}$ be a uniquely split combinatorial noetherian proto-exact category with duality. Then there is a weak homotopy equivalence
\[
\mathcal{GW}(\mathcal{A})
\simeq
\bigsqcup_{w \in W_0(\mathcal{A})}BG_{S_w}\times\mathcal{GW}_H(\mathcal{A}),
\]
where $S_w$ is an isotropically simple representative of the Witt class $w \in W_0(\mathcal{A})$ and $G_{S_w}$ is the self-isometry group of $S_w$.
\end{thm}

See Example \ref{ex: GW-theory of Vect F1} for a description of the spaces $\mathcal{GW}_H(\Vect_\Fun)$, $\mathcal{GW}^{\oplus}(\Vect_\Fun)$ and $\mathcal{GW}(\Vect_\Fun)$.

\subsubsection*{Acknowledgements}
The authors thank Marco Schlichting for helpful correspondence. All three authors thank the Max Planck Institute for Mathematics in Bonn for its hospitality and financial support.


\section{Proto-exact categories with duality}
\label{sec:prelim}

In this section we establish the categorical background required for the remainder of the paper.

\subsection{Proto-exact categories}
\label{sec:cat}

Proto-exact categories, introduced by Dyckerhoff and Kapranov \cite[\S 2.4]{dyckerhoff2019}, are one of many possible non-additive generalizations of Quillen's exact categories \cite[\S 2]{quillen1973}. The axioms of a proto-exact category, which we recall below, are manifestly self-dual and so are particularly well-adapted to Grothendieck--Witt theory.

A proto-exact category is a category $\mathcal{A}$ with a zero object $0$ together with two distinguished classes of morphisms, $\mathfrak{I}$ and $\mathfrak{D}$, called inflations (or admissible monomorphisms) and deflations (or admissible epimorphisms) and denoted $\rightarrowtail$ and $\twoheadrightarrow$, respectively, such that the following axioms hold:
\begin{enumerate}[wide,labelwidth=!, labelindent=0pt,label=(\roman*)]
\item Any morphism $0 \rightarrow U$ is in $\mathfrak{I}$ and any morphism $U \rightarrow 0$ is in $\mathfrak{D}$.

\item The class $\mathfrak{I}$ is closed under composition and contains all isomorphisms, and similarly for $\mathfrak{D}$.

\item A commutative square of the form
\begin{equation}
\label{eq:bicartDiag}
\begin{tikzpicture}[baseline= (a).base]
\node[scale=1] (a) at (0,0){
\begin{tikzcd}
U \arrow[two heads]{d} \arrow[tail]{r} & V \arrow[two heads]{d}\\
W \arrow[tail]{r} & X
\end{tikzcd}
};
\end{tikzpicture}
\end{equation}
is cartesian if and only if it is cocartesian.

\item A diagram of the form $W \rightarrowtail X \twoheadleftarrow V$ can be completed to a bicartesian square of the form \eqref{eq:bicartDiag}.

\item A diagram of the form $W \twoheadleftarrow U \rightarrowtail V$ can be completed to a bicartesian square of the form \eqref{eq:bicartDiag}.
\end{enumerate}

A bicartesian square \eqref{eq:bicartDiag} with $W =0$ is called a conflation and, for ease of notation, is denoted $U \rightarrowtail V \twoheadrightarrow X$. Being kernels, inflations are necessarily monomorphisms. Similarly, deflations are epimorphisms.

If $\mathcal{A}$ is a proto-exact category, then the opposite category $\mathcal{A}^{\op}$ has a natural proto-exact structure in which the inflations (resp. deflations) are the deflations (resp. inflations) in $\mathcal{A}$. The cartesian product of proto-exact categories has a natural proto-exact structure. 

A functor between proto-exact categories is called proto-exact if it preserves zero objects and bicartesian squares of the form \eqref{eq:bicartDiag}. In particular, a proto-exact functor sends conflations to conflations.

\begin{Def}
An exact direct sum on a proto-exact category $\mathcal{A}$ is a symmetric monoi\-dal structure $\oplus$ on $\mathcal{A}$ subject to the following axioms:
\begin{enumerate}[wide,labelwidth=!, labelindent=0pt,label=(DS\arabic*)]
\item \label{ax:monUnit} The monoidal unit of $\mathcal{A}$ is $0$.

\item \label{ax:exFun} The bifunctor $\oplus: \mathcal{A} \times \mathcal{A} \rightarrow \mathcal{A}$ is proto-exact.

\vspace{10pt}
Given objects $U,V \in \mathcal{A}$, set $i_U: U \xrightarrow[]{\id_U \oplus 0_{0 \rightarrowtail V}} U \oplus V$ and $\pi_U: U \oplus V \xrightarrow[]{\id_U \oplus 0_{V \twoheadrightarrow 0}} U$, where we have used \ref{ax:monUnit} to identify $U \oplus 0 \simeq U$ and $0 \oplus V \simeq V$. Axiom \ref{ax:exFun} implies that $i_U$ is an inflation and $\pi_U$ is a deflation.
\vspace{10pt}

\item \label{ax:restrBij} The maps
\[
\Hom_{\mathcal{A}}(U \oplus V, W) \rightarrow \Hom_{\mathcal{A}}(U,W) \times \Hom_{\mathcal{A}}(V,W),
\qquad
f \mapsto (f \circ i_U, f \circ i_V)
\]
and
\[
\Hom_{\mathcal{A}}(W, U \oplus V) \rightarrow \Hom_{\mathcal{A}}(W,U) \times \Hom_{\mathcal{A}}(W,V),
\qquad
f \mapsto (\pi_U \circ f, \pi_V \circ f)
\]
are injections for all $U,V,W \in \mathcal{A}$.

\item \label{ax:sectSplit} Let $U \yrightarrowtail[]{i} X \ytwoheadrightarrow[]{\pi} V$ be a conflation. For each section $s$ of $\pi$, there exists a unique isomorphism $\phi$ which makes the following diagram commute:
\[
\begin{tikzpicture}[baseline= (a).base]
\node[scale=1] (a) at (0,0){
\begin{tikzcd}
& X & \\
U \arrow[tail]{r}[below]{i_U} \arrow[tail]{ur}[above]{i}& U \oplus V \arrow{u}[left]{\phi}& \arrow[tail]{l}[below]{i_V} \arrow{ul}[above]{s} V.
\end{tikzcd}
};
\end{tikzpicture}
\]
For each retraction $r$ of $i$, there exists a unique isomorphism $\psi$ which makes the following diagram commute:
\[
\begin{tikzpicture}[baseline= (a).base]
\node[scale=1] (a) at (0,0){
\begin{tikzcd}
& X \arrow{dl}[above]{r} \arrow[two heads]{dr}[above]{\pi} & \\
U & U \oplus V \arrow[two heads]{r}[below]{\pi_V} \arrow[two heads]{l}[below]{\pi_U} \arrow{u}[left]{\psi}& V.
\end{tikzcd}
};
\end{tikzpicture}
\]
\end{enumerate}
\end{Def}

\begin{Rems}
\label{rem:noUniversal}
\begin{enumerate}[wide,labelwidth=!, labelindent=0pt,label=(\roman*)]
\item Motivated by the theory of Hall algebras, the notion of a proto-exact category with exact direct sum was introduced in \cite[\S 1.1]{mbyoung2021}, where only axioms \ref{ax:monUnit} and \ref{ax:exFun} were imposed. As all categories studied in \cite{mbyoung2021} also satisfy axioms \ref{ax:restrBij} and \ref{ax:sectSplit}, we do not introduce new terminology.

\item Axiom \ref{ax:restrBij} allows for the use of matrix notation to describe morphisms between finite direct sums. Namely, associated to a morphism $f: \bigoplus_i X_i\to\bigoplus_j Y_j$ is the matrix $(\pi_{Y_j} \circ f \circ i_{X_i})_{i,j}$ and this matrix uniquely determines $f.$ However, not every matrix of morphisms in $\mathcal{A}$ arises from a morphism in $\mathcal{A}$. To construct a morphism between direct sums, one can use the symmetric monoidal structure, which allows to construct morphisms from generalized monomial matrices whose entries are morphisms in $\mathcal{A}$, or axiom \ref{ax:sectSplit}.

\item Closely related to proto-exact categories with exact direct sum are Deitmar's quasi-exact categories \cite[\S 3.2]{deitmar2006} and belian categories \cite{deitmar2012}. A key difference between proto-exact categories with exact direct sum and quasi-exact or belian categories is that the exact direct sum $\oplus$ is neither required to be a product nor a coproduct. For this reason, the opposite of a quasi-exact category need not be quasi-exact, and similarly for belian categories, whereas if $\mathcal{A}$ is proto-exact with exact direct sum $\oplus$, then $\mathcal{A}^{\op}$ is proto-exact with exact direct sum $\oplus^{\op}$. Axiom \ref{ax:sectSplit} can be seen as a weak replacement of the universal property of a product and coproduct.
\end{enumerate}
\end{Rems}

A functor between proto-exact categories with exact direct sum is called exact if it is proto-exact and $\oplus$-monoidal. For ease of notation, we do not introduce notation for the $\oplus$-monoidal data of an exact functor.

In the remainder of this section, we prove some basic results about proto-exact categories with exact direct sum.

Let $\mathcal{A}$ be a proto-exact category with exact direct sum. A commutative diagram
\[
\begin{tikzpicture}[baseline= (a).base]
\node[scale=1] (a) at (0,0){
\begin{tikzcd}
U \arrow[tail]{r}[above]{i} \arrow[d,equal] & X \arrow[two heads]{r}[above]{\pi} & V \arrow[d,equal]\\
U \arrow[tail]{r}[below]{i_U} & U \oplus V \arrow{u}[left]{\phi} \arrow[two heads]{r}[below]{\pi_V}& V
\end{tikzcd}
};
\end{tikzpicture}
\]
of conflations with $\phi$ an isomorphism is called a splitting of $U \yrightarrowtail[]{i} X \ytwoheadrightarrow[]{\pi} V$.

\begin{Lem}
\label{lem:uniqueSect}
Let $\mathcal{A}$ be a proto-exact category with exact direct sum. Then there is a bijection between the set of splittings of a conflation $U \yrightarrowtail[]{i} X \ytwoheadrightarrow[]{\pi} V$, the set of sections of $\pi$ and the set of retractions of $i$.
\end{Lem}

\begin{proof}
Given a splitting $\phi$, the map $\phi \circ i_V$ is a section of $\pi$, as follows from the definitions of $\pi_V$ and $i_V$. Given a section $s$ of $\pi$, let $\phi : U \oplus V \rightarrow X$ be the isomorphism whose existence is guaranteed by axiom \ref{ax:sectSplit}. Since $\pi \circ \phi \circ i_U = \pi \circ i =0$ and $\pi \circ \phi \circ i_V = \pi \circ s =\id_V$, axiom \ref{ax:restrBij} yields $\pi \circ \phi = \pi_V$. In particular, $\phi$ is a splitting. That these constructions are mutually inverse follows from the uniqueness in axiom \ref{ax:sectSplit}.

The statements involving retractions can be proven similarly.
\end{proof}

\begin{Lem}
\label{lem:dirSumRestr}
Let $\phi_i : U_i \rightarrow V_i$, $i=1,2$, be morphisms in a proto-exact category with exact direct sum. 
\begin{enumerate}[wide,labelwidth=!, labelindent=0pt,label=(\roman*)]
\item \label{lem:dirSumRestrFirst} The square
\[
\begin{tikzpicture}[baseline= (a).base]
\node[scale=1] (a) at (0,0){
\begin{tikzcd}
U_1 \arrow[tail]{r}[above]{i_{U_1}} \arrow{d}[left]{\phi_1} & U_1 \oplus U_2 \arrow{d}[right]{\phi_1 \oplus \phi_2}\\
V_1 \arrow[tail]{r}[below]{i_{V_1}} & V_1 \oplus V_2
\end{tikzcd}
};
\end{tikzpicture}
\]
commutes, as does the analogous square involving the deflations $\pi_{(-)}$.

\item The morphism $\phi_1 \oplus \phi_2$ is an isomorphism if and only if $\phi_1$ and $\phi_2$ are isomorphisms. 
\end{enumerate}
\end{Lem}

\begin{proof}
The square of the first statement is the outside of the diagram
\[
\begin{tikzpicture}[baseline= (a).base]
\node[scale=1] (a) at (0,0){
\begin{tikzcd}
U_1 \arrow[tail]{r}[above]{i_{U_1}} \arrow{d}[left]{\phi_1} & U_1 \oplus U_2 \arrow{d}[right]{\phi_1 \oplus \id_{U_2}}\\
V_1 \arrow[tail]{r}[below]{i_{V_1}} \arrow[d,equal] & V_1 \oplus U_2 \arrow{d}[right]{\id_{V_1} \oplus \phi_2}\\
V_1 \arrow[tail]{r}[below]{i_{V_1}} & V_1 \oplus V_2.
\end{tikzcd}
};
\end{tikzpicture}
\]
The top and bottom squares commute by the definitions of the inflations $i_{(-)}$.

For the second statement, let $\phi_1 \oplus \phi_2$ be an isomorphism with inverse $f$, so that $\id_{U_1 \oplus U_2} = f \circ (\phi_1 \oplus \phi_2)$ and $\id_{V_1 \oplus V_2}=(\phi_1 \oplus \phi_2) \circ f$. These equations, together with the first part of the lemma, give (using that $\pi_{U_i} \circ i_{U_i} = \id_{U_i}$)
\[
\id_{U_i}
=
\pi_{U_i} \circ f \circ (\phi_1 \oplus \phi_2) \circ i_{U_i}
=
\pi_{U_i} \circ f \circ i_{V_i} \circ \phi_i
\]
and
\[
\id_{V_i}
=
\pi_{V_i} \circ (\phi_1 \oplus \phi_2) \circ f \circ i_{V_i}
=
\phi_i \circ \pi_{U_i} \circ f \circ i_{V_i}.
\]
This shows that $\phi_i$ is an isomorphism with inverse $\pi_{U_i} \circ f \circ i_{V_i}$. The other direction is clear.
\end{proof}

\begin{Lem}
\label{lem:pullbackDirectSum}
Let $\mathcal{A}$ be a proto-exact category with exact direct sum. For any inflation $j: U \rightarrowtail V$, the diagram
\[
\begin{tikzpicture}[baseline= (a).base]
\node[scale=1] (a) at (0,0){
\begin{tikzcd}[column sep=50pt]
U \oplus W \arrow[two heads]{d}[left]{\pi_U} \arrow[tail]{r}[above]{j \oplus \id_W} & V \oplus W \arrow[two heads]{d}[right]{\pi_V}\\
U \arrow[tail]{r}[below]{j} & V
\end{tikzcd}
};
\end{tikzpicture}
\]
is bicartesian. In particular, there is an isomorphism $U \oplus W \simeq U \times_V (V \oplus W)$.
\end{Lem}

\begin{proof}
The diagram in question commutes by Lemma \ref{lem:dirSumRestr}\ref{lem:dirSumRestrFirst}. By the axioms of a proto-exact category, it suffices to prove that the diagram is cocartesian. Consider a commutative diagram
\[
\begin{tikzpicture}[baseline= (a).base]
\node[scale=1] (a) at (0,0){
\begin{tikzcd}[column sep=40pt,row sep=20pt]
U \oplus W \arrow[two heads]{d}[left]{\pi_U} \arrow[tail]{r}[above]{j \oplus \id_W} & V \oplus W \arrow[two heads]{d}[right]{\pi_V} \ar[ddr, "r" above right, bend left=15] & {} \\
U \ar[drr, "l" below left, bend left=-15] \arrow[tail]{r}[below]{j} & V & {} \\
{} & {} & T.
\end{tikzcd}
};
\end{tikzpicture}
\]
Since $r \circ (j \oplus \id_W) = l \circ \pi_U$, we see that $0 = r \circ (j \oplus \id_W) \circ i_W = r \circ i_W$. Hence, axiom \ref{ax:restrBij} implies that $r$ is determined by $r \circ i_V$ through $r = r \circ i_V \circ \pi_V$. We claim that $u : V \xrightarrow[]{r \circ i_V} T$ exhibits the universal property of a cocartesian diagram. We have $u \circ \pi_V = r \circ i_V \circ \pi_V = r$ and
\begin{eqnarray*}
u \circ j \circ \pi_U
&=&
u \circ \pi_V \circ (j \oplus \id_W) \\
&=&
r \circ i_V \circ \pi_V \circ (j \oplus \id_W) \\
&=&
r \circ (j \oplus \id_W) \\
&=&
l \circ \pi_U.
\end{eqnarray*}
Since $\pi_U$ is an epimorphism, we conclude that $u \circ j = l$.
\end{proof}

\begin{Def}
A proto-exact category with exact direct sum is called split (resp. uniquely split) if every conflation admits a splitting (resp. unique splitting).
\end{Def}

A non-zero exact category $\mathcal{A}$ is never uniquely split. Indeed, for a non-zero object $U \in \mathcal{A}$, the set of splittings of the split conflation $U \rightarrowtail U^{\oplus 2} \twoheadrightarrow U$ is a torsor for the additive group $\End_{\mathcal{A}}(U)$. The uniquely split property is therefore only of interest for non-exact proto-exact categories.

\begin{Lem}
\label{lem:uniqueMorphismConflation}
Let $\mathcal{A}$ be a uniquely split proto-exact category. Given a diagram
\[
\begin{tikzpicture}[baseline= (a).base]
\node[scale=1] (a) at (0,0){
\begin{tikzcd}
U_1 \arrow[tail]{r}[above]{i_1} \arrow{d}[left]{f} & E_1 \ar[d,dashed,"h"] \arrow[two heads]{r}[above]{\pi_1}  & V_1 \arrow{d}[right]{g}\\
U_2 \arrow[tail]{r}[below]{i_2} & E_2 \arrow[two heads]{r}[below]{\pi_2} & V_2,
\end{tikzcd}
};
\end{tikzpicture}
\]
in $\mathcal{A}$ in which both rows are conflations, there exists a unique morphism $h:E_1 \rightarrow E_2$ which makes the diagram commute.
\end{Lem}

\begin{proof}
By the uniquely split assumption, we can uniquely extend the given diagram to a commutative diagram
\[
\begin{tikzpicture}[baseline= (a).base]
\node[scale=1.0] (a) at (0,0){
\begin{tikzcd}[column sep=30pt, row sep=30pt]
U_1 \arrow[tail]{r}[above]{i_{U_1}} \arrow[d,equal] & U_1 \oplus V_1 \arrow[two heads]{r}[above]{\pi_{V_1}} \arrow{d}[left]{\phi_1} & V_1 \arrow[d,equal]\\
U_1 \arrow[tail]{r}[above]{i_1} \arrow{d}[left]{f} & E_1 \arrow[two heads]{r}[above]{\pi_1} & V_1 \arrow{d}[right]{g} \\
U_2 \arrow[tail]{r}[above]{i_2} \arrow[d,equal] & E_2 \arrow[two heads]{r}[above]{\pi_2} \arrow{d}[left]{\phi_2} & V_2 \arrow[d,equal] \\
U_2 \arrow[tail]{r}[below]{i_{U_2}} & U_2 \oplus V_2 \arrow[two heads]{r}[below]{\pi_{V_2}} & V_2
\end{tikzcd}
};
\end{tikzpicture}
\]
for some isomorphisms $\phi_1$ and $\phi_2$. The reduces the problem to the case
\[
\begin{tikzpicture}[baseline= (a).base]
\node[scale=1] (a) at (0,0){
\begin{tikzcd}
U_1 \arrow[tail]{r}[above]{i_{U_1}} \arrow{d}[left]{f} & U_1 \oplus V_1 \arrow[two heads]{r}[above]{\pi_{V_1}}  & V_1 \arrow{d}[right]{g}\\
U_2 \arrow[tail]{r}[below]{i_{U_2}} & U_2 \oplus V_2 \arrow[two heads]{r}[below]{\pi_{V_2}} & V_2,
\end{tikzcd}
};
\end{tikzpicture}
\]
in which we can take $h = f \oplus g$. That $h$ is the unique admissible choice follows from axiom \ref{ax:restrBij}.
\end{proof}

In particular, it follows from Lemma \ref{lem:uniqueMorphismConflation} that for a diagram of conflations of the form
\[
\begin{tikzpicture}[baseline= (a).base]
\node[scale=1] (a) at (0,0){
\begin{tikzcd}
U \arrow[tail]{r}[above]{i_1} \arrow[d,equal] & E_1 \arrow[two heads]{r}[above]{\pi_1}  & V \arrow[d,equal]\\
U \arrow[tail]{r}[below]{i_2} & E_2 \arrow[two heads]{r}[below]{\pi_2} & V,
\end{tikzcd}
};
\end{tikzpicture}
\]
there exists a unique morphism $h: E_1 \rightarrow E_2$ making the diagram commute and that, moreover, that $h$ is an isomorphism.

Motivated by \cite{eppolito2020}, we make the following definition. See also \cite[\S 1.1]{mbyoung2021}.

\begin{Def}
A proto-exact category with exact direct sum is called combinatorial if, for each inflation $i: U \rightarrowtail X_1 \oplus X_2$, there exist inflations $i_k: U_k \rightarrowtail X_k$, $k=1,2$, and an isomorphism $f: U \rightarrow U_1 \oplus U_2$ such that $i = (i_1 \oplus i_2) \circ f$ and, dually, for each deflation $\pi: X_1 \oplus X_2 \twoheadrightarrow U$, there exist deflations $\pi_k: X_k \rightarrowtail U_k$, $k=1,2$, and an isomorphism $g: U_1 \oplus U_2 \rightarrow U$ such that $\pi = g \circ (\pi_1 \oplus \pi_2)$.
\end{Def}

Given an inflation $i: U \rightarrowtail X_1 \oplus X_2$, we sometimes write $U \cap X_k$ for the object $U_k$, $k=1,2$.

\begin{Lem}
\label{lem:partialIsoComb}
Let $\mathcal{A}$ be uniquely split combinatorial proto-exact category. If $\phi: U \rightarrowtail X \oplus Y$ is an inflation such that $\pi_X \circ \phi$ is an inflation, then $\pi_Y \circ \phi =0$. Dually, if $\phi: X \oplus Y \twoheadrightarrow U$ is a deflation such that $\phi \circ i_X$ is a deflation, then $\phi \circ i_Y =0$.
\end{Lem}

\begin{proof}
The combinatorial assumption implies that there is a decomposition $U \simeq U_X \oplus U_Y$ under which $\phi$ becomes a morphism
\[
\phi_X \oplus \phi_Y :
U_X \oplus U_Y \rightarrowtail X \oplus Y.
\]
Then we have $\pi_X \circ \phi = \phi_X \circ \pi_{U_X}$. Since $\pi_X \circ \phi$ is an inflation, the kernel $U_Y$ of $\pi_{U_X} : U_X \oplus U_Y \twoheadrightarrow U_X$ is trivial. It follows that $\pi_Y \circ \phi = \phi_Y \circ \pi_{U_Y}=0$. The second statement can be proved in the same way.
\end{proof}

\begin{ex}\label{ex: Vect F1 as proto-exact category}
In our companion paper \cite{eberhardtLorscheidYoung2020b}, we study proto-exact categories occurring in $\Fun$-geometry. These categories come typically with an exact direct sum and are uniquely split and combinatorial. We discuss this in the simplest case of a $\Fun$-linear category, which is the category $\Vect_\Fun$ of ``$\Fun$-vector spaces''. Its objects are pointed sets, which come with a tautological action of the pointed monoid $\Fun=\{0,1\}$. Its morphisms are base point preserving maps that are injective outside the fibre over the base point.
 
 The direct sum $N\oplus M=N\vee M$ of two pointed sets $N$ and $M$ is the wedge sum, that is, their disjoint union modulo the identification the base points. It comes with canonical inclusions and projections, the latter contracting the opposite summand to the base point. It is an easy exercise to verify that every conflation in $\Vect_\Fun$ is of the form $$N\xrightarrowtail[]{i_N} M\oplus N\xtwoheadrightarrow[]{\pi_M} M,$$ so that that $\Vect_\Fun$ is uniquely split. Similarly, it is easy to see that $\Vect_\Fun$ is combinatorial.
\end{ex}

\subsection{Proto-exact categories with duality}
\label{sec:catWD}

For a detailed introduction to categories with duality, the reader is referred to \cite[\S 2]{schlichting2010}.

A category with duality is a triple $(\mathcal{A},P, \Theta)$ consisting of a category $\mathcal{A}$, a functor $P: \mathcal{A}^{\op} \rightarrow \mathcal{A}$ and a natural isomorphism $\Theta: \id_{\mathcal{A}} \Rightarrow P \circ P^{\op}$ which satisfies
\[
P(\Theta_U) \circ \Theta_{P(U)} = \id_{P(U)},
\qquad
U \in \mathcal{A}.
\]
We often omit the pair $(P,\Theta)$ from the notation if it will not cause confusion. If $\mathcal{A}$ has a proto-exact structure and $P$ is proto-exact, then we call this a proto-exact category with duality.

A (nondegenerate) symmetric form in $(\mathcal{A}, P, \Theta)$ is a pair $(M,\psi_M)$ consisting of an object $M \in \mathcal{A}$ and an isomorphism $\psi_M: M \rightarrow P(M)$ which satisfies $P(\psi_M) \circ \Theta_M = \psi_M$. We often write $M$ or $\psi_M$ for $(M,\psi_M)$. An isometry $\phi: (M,\psi_M) \rightarrow (N, \psi_N)$ is an isomorphism $\phi: M \rightarrow N$ which satisfies $\psi_M = P(\phi) \circ \psi_N \circ \phi$. The groupoid of symmetric forms and their isometries is denoted $\mathcal{S}_h$ and called the hermitian groupoid of $\mathcal{A}$.

Let $(M,\psi_M)$  be a symmetric form in a proto-exact category with duality. An inflation $i: U \rightarrowtail M$ is called isotropic if $P(i) \circ \psi_M \circ i$ is zero and the induced monomorphism $U \rightarrow U^{\perp} := \ker(P(i) \circ \psi_M)$ is an inflation. In this case, $M \git U: = U^{\perp} \slash U$, the isotropic reduction of $M$ by $U$, inherits a symmetric morphism $\psi_{M \git U} : M \git U \rightarrow P(M \git U)$.

We say that $(\mathcal{A}, P, \Theta)$ satisfies the Reduction Assumption (as introduced in \cite[\S 3.4]{mbyoung2018b}) if for every symmetric form $(M,\psi_M)$ and every isotropic inflation $i: U \rightarrowtail M$, the symmetric morphism $\psi_{M \git U} : M \git U \rightarrow P(M \git U)$ is an isomorphism. Exact categories satisfy the Reduction Assumption \cite[Lemma 2.6]{schlichting2010}, as do many proto-exact categories \cite{mbyoung2021}, \cite{eberhardtLorscheidYoung2020b}. A symmetric form $(M, \psi_M)$ is called metabolic if it has a Lagrangian, that is, an isotropic subobject $U \rightarrowtail M$ with $U = U^{\perp}$, and is called isotropically simple if it has no non-zero isotropic subobjects.

In the following, we assume that $\mathcal{A}$ has an exact direct sum, that $P$ is exact and that $\Theta$ is a $\oplus$-monoidal natural isomorphism. 
In this case, we have $P(i_U) =\pi_{P(U)}$ and $P(\pi_U) = i_{P(U)}$ for each $U \in \mathcal{A}$. With these assumptions, $\oplus$ defines an orthogonal direct sum of symmetric forms by
\[
(M, \psi_M) \oplus (N, \psi_N)
=
(M \oplus N, \psi_M \oplus \psi_N).
\]
This gives $\mathcal{S}_h$ the structure of a symmetric monoidal groupoid. Moreover, $\oplus$ allows to define hyperbolic symmetric forms. Namely, given an object $U \in \mathcal{A}$, the pair
\[
\left(
H(U) = U \oplus P(U), \psi_{H(U)} = \begin{psmallmatrix} 0 & \id_{P(U)} \\ \Theta_U & 0 \end{psmallmatrix}
\right)
\]
is a symmetric form in $\mathcal{A}$, called the hyperbolic form on $U$. The assignment $U \mapsto H(U)$ extends to a functor $H: \mathcal{S} \rightarrow \mathcal{S}_h$, where $\mathcal{S}$ is the maximal groupoid of $\mathcal{A}$. A symmetric form is called hyperbolic if it is isometric to one of the form $\psi_{H(U)}$.

We now prove an analogue of the fact that split metabolics in an exact category with duality in which ``$2$ is invertible'' are hyperbolic \cite[Lemma 2.5(b)]{hornbostel2002}.

\begin{Lem}
\label{lem:noMeta}
A metabolic form in a uniquely split proto-exact category with duality is hyperbolic.
\end{Lem}

\begin{proof}
Let $i: U \rightarrowtail (M, \psi_M)$ be a Lagrangian. We can extend $i$ to a commutative diagram
\[
\begin{tikzpicture}[baseline= (a).base]
\node[scale=1] (a) at (0,0){
\begin{tikzcd}[column sep=50pt]
U \arrow[tail]{r}[above]{i} \arrow{d}[left]{\Theta_U} & M \arrow[two heads]{r}[above]{P(i) \circ \psi_M} \arrow{d}[left]{\psi_M}  & P(U) \arrow[d,equal]\\
P^2(U) \arrow[tail]{r}[below]{\psi_M \circ i \circ \Theta_U^{-1}} & P(M) \arrow[two heads]{r}[below]{P(i)} & P(U).
\end{tikzcd}
};
\end{tikzpicture}
\]
Since $\mathcal{A}$ is (uniquely) split, this diagram can be extended to the commutative diagram
\[
\begin{tikzpicture}[baseline= (a).base]
\node[scale=1.0] (a) at (0,0){
\begin{tikzcd}[column sep=35pt, row sep=20pt]
U \arrow[tail]{r}[above]{i_U} \arrow[d,equal] & U \oplus P(U) \arrow[two heads]{r}[above]{\pi_U} \arrow{d}[left]{\phi}  & P(U) \arrow[d,equal]\\
U \arrow[tail]{r}[above]{i} \arrow{d}[left]{\Theta_U} & N \arrow[two heads]{r}[above]{P(i) \circ \psi_M} \arrow{d}[left]{\psi_M}  & P(U) \arrow[d,equal]\\
P^2(U) \arrow[tail]{r}[above]{\psi_M \circ i \circ \Theta_U^{-1}} \arrow[d,equal] & P(M) \arrow[two heads]{r}[above]{P(i)} \arrow{d}[left]{P(\phi)} & P(U) \arrow[d,equal] \\
P^2(U) \arrow[tail]{r}[below]{P(\pi_U)}& P(U \oplus P(U)) \arrow[two heads]{r}[below]{P(i_U)} & P(U)
\end{tikzcd}
};
\end{tikzpicture}
\]
for some isomorphism $\phi$. Then $P(\phi) \circ \psi_M \circ \phi$ and $\psi_{H(U)}$ are endomorphisms of the split conflation on $U$ and $P(U)$ and so, by Lemma \ref{lem:uniqueMorphismConflation}, are equal. Hence, $\phi$ is an isometry $\psi_{H(U)} \xrightarrow[]{\sim} \psi_M$.
\end{proof}

\begin{Prop}
\label{prop:splitOffNondegen}
Let $(N, \psi_N)$ be a symmetric form in a uniquely split proto-exact category with duality and $i: M \rightarrowtail N$ an inflation such that $\psi_M := P(i) \circ \psi_N \circ i$ is a symmetric form on $M$. Then there exists a symmetric form $(M^{\prime},\psi_{M^{\prime}})$ and an isometry $(N,\psi_N) \simeq (M,\psi_M) \oplus (M^{\prime},\psi_{M^{\prime}})$ which identifies $i$ with $i_M: M \rightarrowtail M \oplus M^{\prime}$.
\end{Prop}

\begin{proof}
Set
\[
\pi : = \psi_M^{-1} \circ P(i) \circ \psi_N : N \twoheadrightarrow M
\]
and define $M^{\prime} \in \mathcal{A}$ by the conflation $M^{\prime} \yrightarrowtail[]{j} N \ytwoheadrightarrow[]{\pi} M$. Set $\psi_{M^{\prime}} = P(j) \circ \psi_N \circ j$. Fix a splitting
\[
\begin{tikzpicture}[baseline= (a).base]
\node[scale=1] (a) at (0,0){
\begin{tikzcd}[column sep=30pt, row sep=30pt]
M^{\prime} \arrow[tail]{r}[above]{j} \arrow[d,equal] & N \arrow[two heads]{r}[above]{\pi} & M \arrow[d,equal]\\
M^{\prime} \arrow[tail]{r}[below]{i_{M^{\prime}}} & M^{\prime} \oplus M \arrow[two heads]{r}[below]{\pi_M} \arrow{u}[left]{\phi} & M
\end{tikzcd}
};
\end{tikzpicture}
\]
and define a symmetric form on $M^{\prime} \oplus M$ by $\psi = P(\phi) \circ \psi_N \circ \phi$. With these definitions, we have $P(i_{M^{\prime}}) \circ \psi \circ i_{M^{\prime}} = \psi_{M^{\prime}}$ and
\[
P(i_M) \circ \psi \circ i_M
=
P(\phi \circ i_M) \circ \psi_N \circ (\phi \circ i_M)
=
\psi_M.
\]
To see the last equality, note that since $\phi \circ i_M$ and $i$ are sections of $\pi$, uniqueness of splittings and Lemma \ref{lem:uniqueSect} combine to give $\phi \circ i_M = i$. We also have
\[
P(i_M) \circ \psi \circ i_{M^{\prime}}
=
P(\phi \circ i_M) \circ \psi_N \circ j
=
P(i) \circ \psi_N \circ j
=
0
\]
and, dually, $P(i_{M^{\prime}}) \circ \psi \circ i_M=0$. Using axiom \ref{ax:restrBij}, we conclude that $\psi = \psi_M \oplus \psi_{M^{\prime}}$. Finally, $\psi_{M^{\prime}}$ is symmetric by construction and an isomorphism by Lemma \ref{lem:dirSumRestr}. Hence, $(M^{\prime}, \psi_{M^{\prime}})$ is a symmetric form.
\end{proof}

\begin{Prop}
\label{prop:isoReductionSplitting}
Let $(N,\psi_N)$ be a symmetric form in a uniquely split proto-exact category with duality. If $i: U \rightarrowtail (N, \psi_N)$ is isotropic, then there exists an isometry $\phi: N \rightarrow H(U) \oplus ( N \git U)$ which identifies
\[
U \rightarrowtail U^{\perp} \rightarrowtail N
\]
with
\[
U \xrightarrowtail[]{i_U} U \oplus (N \git U) \xrightarrowtail[]{i_{U \oplus (N \git U)}} H(U) \oplus ( N \git U).
\]
\end{Prop}

\begin{proof}
Set $M = N \git U$. Complete $i$ to the following diagram, all of whose rectangles are bicartesian:
\[
\begin{tikzpicture}[baseline= (a).base]
\node[scale=1.0] (a) at (0,0){
\begin{tikzcd}[column sep=35pt, row sep=35pt]
U \arrow[tail]{r}[above]{k} \arrow[two heads]{d} & U^{\perp} \arrow[tail]{r}[above]{j} \arrow[two heads]{d}[left]{\pi} & N \arrow[two heads]{d}[right]{P(j) \circ \psi_N}\\
0 \arrow[tail]{r} & M \arrow[tail]{r}[below]{P(\pi) \circ \psi_M} \arrow[two heads]{d} & P(U^{\perp}) \arrow[two heads]{d}[right]{P(k)} \\
{} & 0 \arrow[tail]{r} & P(U).
\end{tikzcd}
};
\end{tikzpicture}
\]
Here $i = j \circ k$. Let $f: U^{\perp} \xrightarrow[]{\sim} U \oplus M$ be the unique splitting of $U \rightarrowtail U^{\perp} \twoheadrightarrow M$ and set $j^{\prime} = j \circ f^{-1}$. Let $\phi: N \xrightarrow[]{\sim} U \oplus P(U) \oplus M$ be the unique splitting of $U \oplus M \xrightarrowtail[]{j^{\prime}} N \twoheadrightarrow P(U)$ and set $\psi^{\prime} = P(\phi^{-1}) \circ \psi_N \circ \phi^{-1}$. Then the previous diagram becomes
\[
\begin{tikzpicture}[baseline= (a).base]
\node[scale=1.0] (a) at (0,0){
\begin{tikzcd}[column sep=35pt, row sep=35pt]
U \arrow[tail]{r}[above]{i_U} \arrow[two heads]{d} & U \oplus M \arrow[tail]{r}[above]{i_{U \oplus M}} \arrow[two heads]{d}[left]{\pi_M} & U \oplus P(U) \oplus M \arrow[two heads]{d}[right]{P(i_{U \oplus M}) \circ \psi^{\prime}}\\
0 \arrow[tail]{r} & M \arrow[tail]{r}[below]{P(\pi_M) \circ \psi_M} \arrow[two heads]{d} & P(U \oplus M) \arrow[two heads]{d}[right]{P(i_U)} \\
{} & 0 \arrow[tail]{r} & P(U).
\end{tikzcd}
};
\end{tikzpicture}
\]
In particular, we have
\[
P(\pi_M) \circ \psi_M \circ \pi_M
=
P(i_{U \oplus M}) \circ \psi^{\prime} \circ i_{U \oplus M}.
\]
Pre- and post-composing this equation with $i_M : M \rightarrowtail U \oplus M$ and $P(i_M)$, respectively, gives
\[
P(i_M) \circ P(\pi_M) \circ \psi_M \circ \pi_M \circ i_M
=
P(i_M) \circ P(i_{U \oplus M}) \circ \psi^{\prime} \circ i_{U \oplus M} \circ i_M,
\]
which can be rewritten as
\[
\psi_M
=
P(i_{M \rightarrowtail U\oplus P(U) \oplus M}) \circ \psi^{\prime} \circ i_{M \rightarrowtail U\oplus P(U) \oplus M}.
\]
Proposition \ref{prop:splitOffNondegen} therefore implies that there exists an isometry
\[
(U \oplus P(U) \oplus M, \psi) \simeq (U \oplus P(U), \psi_{U \oplus P(U)}) \oplus (M, \psi_M)
\]
under which $U \rightarrowtail U \oplus P(U) \oplus M$ factors through the standard Lagrangian $i_U: U \rightarrowtail (U \oplus P(U), \psi_{U \oplus P(U)})$. We conclude using Lemma \ref{lem:noMeta} that $\psi_{U \oplus P(U)} \simeq \psi_{H(U)}$.
\end{proof}

A proto-exact category is called noetherian if any ascending chain of inflations stabilizes after finitely many steps.

\begin{Prop}
\label{prop:isoSimpDecomp}
Let $(M, \psi_M)$ be a symmetric form in a uniquely split noetherian proto-exact category $\mathcal{A}$ with duality. 
\begin{enumerate}[wide,labelwidth=!, labelindent=0pt,label=(\roman*)]
\item \label{propPart:essSurj} There exists an object $U \in \mathcal{A}$, an isotropically simple symmetric form $(N, \psi_N) \in \mathcal{S}_h$ and an isometry $M \simeq H(U) \oplus N$.

\item \label{propPart:fullFaith} If, moreover, $\mathcal{A}$ is combinatorial,
then the decomposition from part (i) is unique up to isometry in either summand.
\end{enumerate}
\end{Prop}

\begin{proof}
If $M$ is isotropically simple, then take $U=0$ and we are done. Otherwise, let $U_1 \rightarrowtail M$ be a non-zero isotropic. Then $M \simeq H(U_1) \oplus M_2$ by Proposition \ref{prop:isoReductionSplitting}. If $M_2$ is isotropically simple, then we are done. Otherwise, choose a non-zero isotropic $\tilde{U}_2 \rightarrowtail M_2$, so that $M \simeq H(U_2) \oplus M_3$ with $U_2 = U_1 \oplus \tilde{U_2}$. Continuing in this way, we obtain an ascending chain of inflations $U_1 \rightarrowtail U_2 \rightarrowtail \cdots$ of $M$, which is finite by the noetherian hypothesis. This proves the first part.

For the second part, let $\phi : H(U) \oplus N \rightarrow H(U^{\prime}) \oplus N^{\prime}$ be an isometry with $N$, $N^{\prime}$ isotropically simple. Then $\phi \circ i_U: U \rightarrowtail H(U^{\prime}) \oplus N^{\prime}$ is isotropic and hence so too are the subobjects $\phi(U) \cap H(U^{\prime})$ and $\phi(U) \cap N^{\prime}$. Since $N^{\prime}$ is isotropically simple, $\phi(U) \cap N^{\prime} =0$, whence $\phi \circ i_U$ factors through $H(U^{\prime})$. Similarly, $\phi \circ i_{P(U)}$ factors through $H(U^{\prime})$. We therefore have a commutative diagram
\begin{equation}
\label{diag:combSplitting}
\begin{tikzpicture}[baseline= (a).base]
\node[scale=1] (a) at (0,0){
\begin{tikzcd}[column sep=30pt, row sep=30pt]
H(U) \arrow[tail]{r}[above]{i_{H(U)}} \arrow{d}[left]{\phi_L} & H(U) \oplus N \arrow[two heads]{r}[above]{\pi_N} \arrow{d}[left]{\phi} & N \arrow{d}[right]{\phi_R}\\
H(U^{\prime}) \arrow[tail]{r}[below]{i_{H(U^{\prime})}} & H(U^{\prime}) \oplus N^{\prime} \arrow[two heads]{r}[below]{\pi_{N^{\prime}}} & N^{\prime}.
\end{tikzcd}
};
\end{tikzpicture}
\end{equation}
By applying the same argument to $\phi^{-1}$ and patching the resulting diagram with diagram \eqref{diag:combSplitting} (similar to the proof of Lemma \ref{lem:noMeta}), we conclude that $\phi_L$ and $\phi_R$ are isomorphisms. Since $\phi$ is symmetric, so too are $\phi_L$ and $\phi_R$.
\end{proof}

\begin{ex}\label{ex: duality for Vect F1}
An exact duality for the category $\Vect_\Fun$ from Example \ref{ex: Vect F1 as proto-exact category} is the functor $P:\Vect^\op_\Fun\to\Vect_\Fun$ that is the identity on objects and maps a morphism $f:M\to N$ to its adjoint $f^\ad:N\to M$, defined by $f^\ad(n)=m$ if $n=f(m)$ and $f^\ad(n)=0$ if $n\notin\im(f)$, where $0$ is the base point of $M$. The natural isomorphism $\Theta$ is the identity between $\id_{\Vect_\Fun}$ and $P\circ P^\op=\id_{\Vect_\Fun}$. 
\end{ex}

\section{Algebraic $K$-theory of uniquely split proto-exact categories}
\label{sec:algKProtoEx}

We study $K$-theory spaces of proto-exact and symmetric monoidal categories. Our main result, Theorem \ref{thm:KComparison}, is a Group Completion Theorem for uniquely split proto-exact categories. The techniques and results of this section are used to study Grothendieck--Witt theory in Section \ref{sec:GWProtoEx}.

\subsection{The proto-exact $Q$-construction}
\label{sec:QConstr}

Let $\mathcal{A}$ be a proto-exact category. Similarly to the exact case \cite[\S 2]{quillen1973}, the Quillen $Q$-construction of $\mathcal{A}$ is the category $Q(\mathcal{A})$ defined as follows. An object of $Q(\mathcal{A})$ is simply an object of $\mathcal{A}$. A morphism $U \rightarrow V$ in $Q(\mathcal{A})$ is an equivalence class of diagrams $U \twoheadleftarrow E \rightarrowtail V$ in $\mathcal{A}$. Two such diagrams $U \twoheadleftarrow E \rightarrowtail V$ and $U \twoheadleftarrow E^{\prime} \rightarrowtail V$ are equivalent if there exists an isomorphism $E \rightarrow E^{\prime}$ in $\mathcal{A}$ which makes the obvious diagrams commute. Composition of morphisms in $Q(\mathcal{A})$ is defined via pullback and is well-defined by the axioms of a proto-exact category. The $K$-theory space of $\mathcal{A}$ is
\[
\mathcal{K}(\mathcal{A}) = \Omega B Q(\mathcal{A}),
\]
the based loop space of the classifying space of $Q(\mathcal{A})$, where $0 \in Q(\mathcal{A})$ determines the basepoint of $B Q(\mathcal{A})$. The algebraic $K$-theory groups of $\mathcal{A}$ are the homotopy groups
\[
K_i(\mathcal{A}) = \pi_i \mathcal{K}(\mathcal{A}),
\qquad
i \geq 0.
\]
We remark that an exact direct sum on $\mathcal{A}$ induces a symmetric monoidal structure on $Q(\mathcal{A})$.

\begin{Ex}
Let $\mathcal{A}$ be a proto-exact category. Then $K_0(\mathcal{A})$ is isomorphic to the free group on isomorphism classes of objects of $\mathcal{A}$ modulo the relation $[V] = [U] [W]$ whenever $U \rightarrowtail V \twoheadrightarrow W$. In general, the group $K_0(\mathcal{A})$ need not be abelian \cite[Example 3.6.3]{hekking2017}. If, however, $\mathcal{A}$ admits an exact direct sum, then by using the split conflation on $U$ and $W$ we see that $[U] [W] = [W] [U]$. In this case, $K_0(\mathcal{A})$ is abelian and is described in the familiar way. 
\end{Ex}

\subsection{Direct sum $K$-theory}
\label{sec:dirSumKThy}

Let $(\mathcal{A}, \oplus)$ be a symmetric monoidal category. Note that proto-exact categories with exact direct sum and exact categories define symmetric monoidal categories by forgetting the (proto-)exact structure. The maximal groupoid $\mathcal{S}$ of $\mathcal{A}$ is symmetric monoidal. The direct sum $K$-theory space of $\mathcal{A}$ is the group completion of $B\mathcal{S}$:
\[
\mathcal{K}^{\oplus}(\mathcal{A})
=
B (\mathcal{S}^{-1} \mathcal{S}).
\]
See, for example, \cite[Page 222]{grayson1976}. We refer the reader to \cite[\S IV.4]{weibel2013} for our conventions on group completion.

\begin{ex}\label{ex: K-theory of Vect F1}
We calculate the $K$-theory of certain $\Fun$-linear categories in our companion paper \cite{eberhardtLorscheidYoung2020b}. In the simplest case, that of $\Vect_\Fun$ (see Example \ref{ex: Vect F1 as proto-exact category}), we find
 \[
  \mathcal{K}^\oplus(\Vect_\Fun) \ \simeq \ \Z\times (B\Sigma_\infty)^+,
 \]
 as a special case of \cite[Theorem 2.5]{eberhardtLorscheidYoung2020b}.
\end{ex}

\subsection{A Group Completion Theorem for $K$-theory of uniquely split proto-exact categories}
\label{sec:Q=+Thm}

Quillen's Group Completion Theorem \cite{grayson1976} states that for any split exact category $\mathcal{A}$, there is a homotopy equivalence $\mathcal{K}(\mathcal{A}) \simeq \mathcal{K}^{\oplus}(\mathcal{A})$. In this section, we prove a proto-exact version of this result.

\begin{Thm}
\label{thm:KComparison}
Let $\mathcal{A}$ be a uniquely split proto-exact category. Then there is a homotopy equivalence $\mathcal{K}(\mathcal{A}) \simeq \mathcal{K}^{\oplus}(\mathcal{A})$.
\end{Thm}

It is claimed (without proof) in \cite[\S 4]{chu2010} that the proof from \cite{grayson1976} for split exact categories can be modified to apply to (not necessarily uniquely) split proto-exact categories. Instead of following this line of thought, we give a simplified proof in the uniquely split case, which is the setting of interest for this paper. These arguments are used in Section \ref{sec:GWProtoEx}.

We begin with some preparatory material. Write $Q$ for the category $Q(\mathcal{A})$ and let $\mathcal{E}=\mathcal{E}(\mathcal{A})$ be the category of conflations in $\mathcal{A}$. Objects of $\mathcal{E}$ are conflations. Morphisms in $\mathcal{E}$ (from primed to unprimed conflations) are equivalence classes of commutative diagrams
\begin{center}
\begin{tikzcd}
A' \arrow[r, tail]                                                 & B' \arrow[r, two heads] \arrow[d, equal] & C'                                        \\
A \arrow[u, tail] \arrow[d, equal] \arrow[r, tail] & B' \arrow[d, tail] \arrow[r, two heads]               & C_1 \arrow[u, two heads] \arrow[d, tail] \\
A \arrow[r, tail]                                                 & B \arrow[r, two heads]                                & C                                       
\end{tikzcd}
\end{center}
in $\mathcal{A}$ whose rows are conflations. The equivalence relation on morphisms is generated by automorphisms of $C_1$. Projection to the right column defines a fibered functor $g: \mathcal{E} \rightarrow Q$. Given $C\in Q$, let $\mathcal{E}_C := g^{-1}(C)$ be the fibre category. The symmetric monoidal groupoid $\mathcal{S}$ acts on $\mathcal{E}$ by
\[
A'\cdot(A\xrightarrowtail[]{i} B \xtwoheadrightarrow[]{\pi} C)=(A'\oplus A \xrightarrowtail[]{\id_{A^{\prime}} \oplus i} A'\oplus B \xtwoheadrightarrow[]{\pi} C).
\]
Since the functor $\oplus$ is proto-exact, the right hand side is indeed a conflation. The $\mathcal{S}$-action restricts to the fibres of $g$ and, when $Q$ is equipped with the trivial $\mathcal{S}$-action, the base change maps of $g$ are $\mathcal{S}$-equivariant. The functor $g$ is therefore cartesian; see \cite[Pages 222 and 226]{grayson1976}.

The uniquely split assumption on $\mathcal{A}$ allows for the following explicit description of the fibers of $g$. We note that there is no such description for exact categories.

\begin{Lem}
\label{lem:descriptionoffibers} The functor $F_C: \mathcal{S}\rightarrow \mathcal{E}_C$, given on objects and morphisms by
\[
A \quad \longmapsto \quad (A\yrightarrowtail[]{i_A} C\oplus A \ytwoheadrightarrow[]{\pi_C} C)
\]
and
\[
(\phi: A^{\prime} \xrightarrow[]{\sim} A)
\quad \longmapsto \quad 
\begin{tikzcd}[column sep=4em]
A' \arrow["i_{A'}",r, tail]                                                 & C\oplus A' \arrow["\pi_{C}",r, two heads] \arrow[d, equal] & C                                       \\
A \arrow["\phi^{-1}",u, tail] \arrow[d, equal] \arrow["i_{A^{\prime}} \circ \phi^{-1}",r, tail] & C\oplus A' \arrow["\id_C \oplus \phi",d, tail] \arrow["\pi_C",r, two heads]               & C \arrow[u, equal] \arrow[d, equal] \\
A \arrow["i_{A}",swap,r, tail]                                                 & C\oplus A  \arrow["\pi_{C}",swap,r, two heads]                                & C                                      
\end{tikzcd}
\]
respectively, is an $\mathcal{S}$-equivariant equivalence.
\end{Lem}

\begin{proof}
Since $\mathcal{A}$ is (uniquely) split, $F_C$ is essentially surjective. That $F_C$ is faithful follows from the definitions. To see that $F_C$ is full, note that a morphism $f: F_C(A')\rightarrow F_C(A)$ in $\mathcal{E}_C$ is represented by a commutative diagram
\begin{center}
\begin{tikzcd}
A' \arrow["i_{A'}",r, tail]                                                 & C\oplus A' \arrow["\pi_{C}",r, two heads] \arrow[d, equal] & C                                       \\
A \arrow["\alpha",u, tail] \arrow[d, equal] \arrow["i",r, tail] & C\oplus A' \arrow["\beta",d, tail] \arrow["\pi_C",r, two heads]               & C \arrow[u, equal] \arrow[d, equal] \\
A \arrow["i_{A}",swap,r, tail]                                                 & C\oplus A  \arrow["\pi_{C}",swap,r, two heads]                                & C.                                       
\end{tikzcd}
\end{center}
Since $i_{A^{\prime}}$ and $i = i_{A^{\prime}} \circ \alpha$ are kernels of $\pi_C : C \oplus A \twoheadrightarrow C$, the map $\alpha$ is an isomorphism. Post-composing the equation $i_A = \beta \circ i_{A^{\prime}} \circ \alpha$ with $\pi_A$ gives $\id_A = (\pi_A \circ \beta \circ i_{A^{\prime}}) \circ \alpha$, from which we conclude that $\alpha^{-1} = \pi_A \circ \beta \circ i_{A^{\prime}}$. It follows that both $\beta$ and $\id_C \oplus \alpha^{-1}$ define morphisms of conflations
\[
\begin{tikzpicture}[baseline= (a).base]
\node[scale=1] (a) at (0,0){
\begin{tikzcd}[column sep=4em]
A^{\prime} \arrow[tail]{r}[above]{i_{A^{\prime}}} \arrow[d,equal] & C \oplus A^{\prime} \arrow[two heads]{r}[above]{\pi_C} \arrow{d}[left]{\beta}[right]{\id_C \oplus \alpha^{-1}}  & C \arrow[d,equal]\\
A^{\prime} \arrow[tail]{r}[below]{i_A \circ \alpha^{-1}} & C \oplus A \arrow[two heads]{r}[below]{\pi_C} & C.
\end{tikzcd}
};
\end{tikzpicture}
\]
Applying Lemma \ref{lem:uniqueMorphismConflation}, we conclude that $\beta = \id_C \oplus \alpha^{-1}$, so that $f = F_C(\alpha^{-1})$. Hence, $F_C$ is an equivalence.

The $\mathcal{S}$-equivariance of $F_C$ follows from the properties of $\oplus$.
\end{proof}

The morphism $z_C = (0 \twoheadleftarrow 0 \rightarrowtail C)$ in $Q$ induces a base change functor $z_C^*: \mathcal{E}_C\rightarrow \mathcal{E}_0$, given on objects by 
\[
z_C^*(A\rightarrowtail B \twoheadrightarrow C) = (A\yrightarrowtail[]{\id_A} A \ytwoheadrightarrow[]{} 0)
\]
and on morphisms in the obvious way. Dually, the morphism $p_C = ( 0\twoheadleftarrow C \yrightarrowtail[]{\id_C} C)$ induces a base change functor $p_C^*: \mathcal{E}_C \rightarrow \mathcal{E}_0$, given on objects by
\[
p_C^*(A\rightarrowtail B \twoheadrightarrow C) = (B\yrightarrowtail[]{\id_B} B\ytwoheadrightarrow[]{} 0).
\]
There is also a functor $p_{C*}:  \mathcal{E}_0\rightarrow \mathcal{E}_C$ given on objects by 
\[
p_{C*}(A\rightarrowtail B \twoheadrightarrow 0)
=
(A\yrightarrowtail[]{} C\oplus B \ytwoheadrightarrow[]{\pi_C} C).
\]

\begin{Lem}\label{lem:propertiesOfBaseChanges}
\phantomsection
\begin{enumerate}[wide,labelwidth=!, labelindent=0pt,label=(\roman*)]
\item The functor $z^*_{C}$ is an equivalence.

\item There are natural isomorphisms of functors
\[
C\cdot \simeq p^*_Cp_{C*}: \mathcal{E}_0 \rightarrow \mathcal{E}_0
\]
and
\[
C\cdot \simeq p_{C*}p^*_C: \mathcal{E}_C \rightarrow \mathcal{E}_C,
\]
where $C \cdot $ denotes the action of $C \in \mathcal{S}$.
\end{enumerate}
\end{Lem}

\begin{proof}
It is straightforward to verify that there is a natural isomorphism $z^*_C F_C \simeq F_0$. Since $F_C$ and $F_0$ are equivalences (Lemma \ref{lem:descriptionoffibers}), so too is $z^*_C$. The second natural isomorphism is proved similarly, using instead that $p^*_Cp_{C*}F_0 \simeq C\cdot F_0$ and $p_{C*}p^*_CF_C \simeq C\cdot F_C$.
\end{proof}

Denote by $\tilde{g}$ the composition $
\mathcal{S}^{-1}\mathcal{E} \xrightarrow[]{\mathcal{S}^{-1}g} \mathcal{S}^{-1}Q \longrightarrow Q$, the second functor being the equivalence induced by the projection $\mathcal{S} \times Q \rightarrow Q$.

\begin{Lem}\label{lem:homotopyFibrationComparisonKTheory}
After passing to classifying spaces, the sequence $\mathcal{S}^{-1} \mathcal{E}_0 \rightarrow \mathcal{S}^{-1} \mathcal{E} \xrightarrow[]{\tilde{g}} Q$

is a homotopy fibration.
\end{Lem}
\begin{proof}
Since $g$ is cartesian, the fibres of $\tilde{g}$ are of the form $\mathcal{S}^{-1}\mathcal{E}_C$ and the base change map for a morphism $f:C'\rightarrow C$ in $Q$ is $\mathcal{S}^{-1}f^*,$ where $f^*:\mathcal{E}_C\rightarrow \mathcal{E}_{C'}$ is the base change for $g$; see \cite[Page 222]{grayson1976}. To show that all base change maps are equivalences it suffices to consider the morphisms $z_C$ and $p_C$; \textit{cf}. \cite[Lemma 7.9]{srinivas2008}. By Lemma \ref{lem:propertiesOfBaseChanges}, the functor $z_C^*$ is an equivalence, hence so too is $\mathcal{S}^{-1}z_C^*$. Moreover, $\mathcal{S}^{-1}p_{C*}$ is quasi-inverse to $\mathcal{S}^{-1}p^*_{C}$ up to multiplication by $C$. Since multiplication by $C$ is an equivalence on the $\mathcal{S}$-localized categories, Quillen's Theorem B \cite[\S 1]{quillen1973} implies that the sequence in question is a homotopy fibration.
\end{proof}

\begin{proof}[Proof of Theorem \ref{thm:KComparison}]
The space $B\mathcal{S}^{-1}\mathcal{E}$ is contractible, as can be seen in the same way as for the exact case \cite[Page 228]{grayson1976}. The homotopy fibration of Lemma \ref{lem:homotopyFibrationComparisonKTheory} therefore implies that $\Omega BQ$ is homotopy equivalent to $B\mathcal{S}^{-1}\mathcal{E}_0$. By Lemma \ref{lem:descriptionoffibers}, the categories $\mathcal{S}^{-1}\mathcal{E}_0$ and $\mathcal{S}^{-1}\mathcal{S}$ are equivalent. The theorem follows.
\end{proof}


\section{Grothendieck--Witt theory of uniquely split proto-exact categories}
\label{sec:GWProtoEx}

The subject of this section is the Grothendieck--Witt spaces of proto-exact and symmetric monoidal categories with duality. Our first result, Theorem \ref{thm:GWComparison}, is a hyperbolic Group Completion Theorem for uniquely split proto-exact categories with duality. We then use this result to describe the connected components of the Grothendieck--Witt space defined via the hermitian $Q$-construction. The result is Theorem \ref{thm:GWComputation}.

\subsection{The proto-exact hermitian $Q$-construction}
\label{sec:hermQConstr}

The Grothendieck--Witt theory of an exact category with duality can be defined using the hermitian $Q$-construction \cite[\S 0]{charney1986}, \cite[\S 3]{uridia1990}; see also \cite[\S 4.1]{schlichting2010}. This definition extends as follows to a proto-exact category with duality $(\mathcal{A},P, \Theta)$ to define a category $Q_h(\mathcal{A})$. An object of $Q_h(\mathcal{A})$ is a symmetric form in $\mathcal{A}$. A morphism $(M, \psi_M) \rightarrow (N, \psi_N)$ in $Q_h(\mathcal{A})$ is an equivalence class of diagrams $M \xtwoheadleftarrow[]{\pi} E \yrightarrowtail[]{j} N$ in $\mathcal{A}$ such that $j$ is coisotropic, that is, $U: = \ker (P(j) \circ \psi_N) \rightarrowtail N$ is isotropic, and the diagram
\[
\begin{tikzpicture}[baseline= (a).base]
\node[scale=1.0] (a) at (0,0){
\begin{tikzcd}[column sep=40pt, row sep=35pt]
E \arrow[tail]{r}[above]{j} \arrow[two heads]{d}[left]{\pi} & N \arrow[two heads]{d}[right]{P(j) \circ \psi_N}\\
M \arrow[tail]{r}[below]{P(\pi) \circ \psi_M} & P(E)
\end{tikzcd}
};
\end{tikzpicture}
\]
is bicartesian. In words, a morphism $M \rightarrow N$ in $Q_h(\mathcal{A})$ is a presentation of $M$ as an isotropic reduction of $N$. Equivalence and composition of morphisms is as in $Q(\mathcal{A})$. Denote by $F: Q_h(\mathcal{A}) \rightarrow Q(\mathcal{A})$ the forgetful functor.

The Grothendieck--Witt space $\mathcal{GW}(\mathcal{A})$ is the homotopy fibre of $BF: B Q_h(\mathcal{A}) \rightarrow B Q(\mathcal{A})$ over $0$. The Grothendieck--Witt groups of $\mathcal{A}$ are the homotopy groups
\[
GW_i(\mathcal{A})= \pi_i \mathcal{GW}(\mathcal{A}),
\qquad
i \geq 0.
\]
Despite the name, we note that, without further assumptions, $GW_0(\mathcal{A})$ is in fact only a pointed set. The higher Witt groups of $\mathcal{A}$ are
\begin{equation}
\label{eq:wittCoker}
W_i(\mathcal{A}) = \coker \left( K_i(\mathcal{A}) \xrightarrow[]{\pi_i H_*} GW_i(\mathcal{A}) \right),
\qquad
i \geq 1
\end{equation}
where $H_* : \mathcal{K}(\mathcal{A})=\Omega B Q(\mathcal{A}) \rightarrow \mathcal{GW}(\mathcal{A})$ is the natural map, defined up to homotopy).

Suppose now that $\mathcal{A}$ has an exact direct sum $\oplus$. This induces symmetric monoidal structures on $\mathcal{S}_h$ and $Q_h(\mathcal{A})$ and $F$ extends to a symmetric monoidal functor. In this case, $GW_0(\mathcal{A})$ is a commutative monoid. Define the Witt monoid $W_0(\mathcal{A})$ to be the monoid $\pi_0(\mathcal{S}_h)$ of isometry classes of symmetric forms modulo the submonoid of metabolic symmetric forms. Then $W_0(\mathcal{A})$ is a commutative monoid. For a comparison of this definition with that of $W_i(\mathcal{A})$, $i \geq 1$, see the comments below Corollary \ref{cor:GW0Description}.

\begin{Ex}
Suppose that $\mathcal{A}$ is an exact category. In this case $W_0(\mathcal{A})$ is a group. Indeed, the inverse of $[(M, \psi_M)] \in W_0(\mathcal{A})$ is $[(M, - \psi_M)]$, since the diagonal map $M \xrightarrowtail[]{\Delta} \psi_M \oplus - \psi_M$ is Lagrangian.  Similarly, $GW_0(\mathcal{A})$ is an abelian group \cite[Proposition 4.11]{schlichting2010}.
\end{Ex}

\begin{Rem}
Alternatively, in Sections \ref{sec:QConstr} and \ref{sec:hermQConstr} we could have used the proto-exact Waldhausen $\mathcal{S}_{\bullet}$-construction $\mathcal{S}_{\bullet}(\mathcal{A})$ \cite[\S 1.3]{waldhausen1985}, \cite[\S 2.4]{dyckerhoff2019} and hermitian $\mathcal{R}_{\bullet}$-construction $\mathcal{R}_{\bullet}(\mathcal{A})$ \cite[\S 1.8]{hornbostel2004}, \cite[\S 3.4] {mbyoung2018b} to define the $K$-theory and Grothendieck--Witt theory spaces of $\mathcal{A}$, respectively. The resulting spaces are homotopy equivalent to those defined using the $Q$-constructions in a way which is compatible with the maps $H_*$ and $F$.
\end{Rem}


%
%

\subsection{Direct sum Grothendieck--Witt theory}
\label{sec:dirSumGWThy}

Let $(\mathcal{A},\oplus)$ be a symmetric monoidal category with duality. Then $\mathcal{S}_h$ is a symmetric monoidal groupoid. Following \cite[\S 2]{hornbostel2002}, the direct sum Grothendieck--Witt space of $\mathcal{A}$ is the group completion
\[
\mathcal{GW}^{\oplus}(\mathcal{A})
=B(\mathcal{S}_h^{-1} \mathcal{S}_h).
\]
The direct sum Grothendieck--Witt and Witt groups are
\[
GW_i^{\oplus}(\mathcal{A}) = \pi_i \mathcal{GW}^{\oplus}(\mathcal{A}),
\qquad
i \geq 0
\]
and
\[
W_i^{\oplus}(\mathcal{A}) = \coker(K^{\oplus}_i(\mathcal{A}) \xrightarrow[]{H} GW_i^{\oplus}(\mathcal{A})),
\qquad
i \geq 0
\]
respectively. More precisely, $H$ is induced by the map $\mathcal{K}^{\oplus}(\mathcal{A}) \rightarrow \mathcal{GW}^{\oplus}(\mathcal{A})$ determined by the symmetric monoidal hyperbolic functor $H : \mathcal{S} \rightarrow \mathcal{S}_h$. Note that $GW^{\oplus}_i(\mathcal{A})$ and $W_i^{\oplus}(\mathcal{A})$, $i \geq 0$, are abelian groups.

\subsection{A hyperbolic Group Completion Theorem for Grothendieck--Witt theory}
\label{sec:Q=+ThmGW}

In this section we initiate the study of the relation between $\mathcal{GW}(\mathcal{A})$ and $\mathcal{GW}^{\oplus}(\mathcal{A})$ for a uniquely split proto-exact category $\mathcal{A}$. In the split exact case, the spaces $\mathcal{GW}(\mathcal{A})$ and $\mathcal{GW}^{\oplus}(\mathcal{A})$ are homotopy equivalent, thereby giving a hermitian analogue of the Group Completion Theorem; see \cite[Theorem 4.2]{schlichting2004}, \cite[Theorem A.1]{schlichting2017} when ``$2$ is invertible'' and \cite[Theorem 6.6]{schlichting2019} in general. While the naive analogue of the Group Completion Theorem for Grothendieck--Witt theory fails in the proto-exact setting (cf.\ Example \ref{ex: GW-theory of Vect F1}), a hyperbolic analogue does indeed hold.

To begin, we introduce some notation. Let $\mathcal{S}_H \subset \mathcal{S}_h$ and $Q_{H}(\mathcal{A}) \subset Q_h(\mathcal{A})$ be the full subcategories of hyperbolic objects. Denote by $\mathcal{GW}_H(\mathcal{A})$ the homotopy fibre of $BF: BQ_H(\mathcal{A}) \rightarrow BQ(\mathcal{A})$ over $0$ and set $\mathcal{GW}^\oplus_H(\mathcal{A})=B(\mathcal{S}_H^{-1}\mathcal{S}_H)$. We often omit $\mathcal{A}$ from the notation so that, for example, $Q_H = Q_H(\mathcal{A})$.

\begin{Thm}
\label{thm:GWComparison}
Let $\mathcal{A}$ be a uniquely split proto-exact category with duality. Then there is a weak homotopy equivalence $\mathcal{GW}_H(\mathcal{A}) \simeq \mathcal{GW}_H^{\oplus}(\mathcal{A})$.
\end{Thm}

The proof of Theorem \ref{thm:GWComparison} occupies the remainder of this section. 

\begin{Rem}
In \cite[Theorem]{charney1986} a Group Completion Theorem for Grothendieck--Witt theory of exact categories in which $2$ is invertible is claimed. However, the proof contains a serious error which does not appear to be fixable; see \cite[Section 2, Remark]{hornbostel2004}. Correct proofs were found only much later \cite{schlichting2004}, \cite{schlichting2017}, \cite{schlichting2019}. However, since these approaches rely on the additive structure of exact categories in an essential way, they do not carry over to the proto-exact setting. Instead, it turns out that, by restricting to proto-exact categories which are uniquely split, we \emph{can} adapt parts of the strategy of \cite{charney1986} while avoiding the critical mistake.
\end{Rem}

Let $\tau: \mathcal{S}_H \rightarrow  Q_{H}$ be the functor which is the identity on objects and sends an isometry $M \xrightarrow[]{\phi} N$ to $M \xtwoheadleftarrow[]{\id_M} M \xrightarrowtail[]{\phi} N$. We begin by studying the (right) comma categories of $\tau$. Fix $M \in Q_H$. The groupoid $\mathcal{S}$ acts on the comma category $M\backslash \tau$ by
\[
V\cdot (N, (M\stackrel{q}{\twoheadleftarrow} E \stackrel{j}{\rightarrowtail} N))
=
\left(
H(V)\oplus N, (M\xtwoheadleftarrow[]{q \circ \pi_E} V \oplus E \xrightarrowtail[]{i_V \oplus j} H(V)\oplus N)
\right).
\]
The following lemma, which relies on the uniquely split assumption and is not true in the exact setting, gives an explicit description of $M\backslash \tau$.

\begin{Lem}
\label{lem:HMEquiv}
The functor $H^M:\mathcal{S} \rightarrow M\backslash\tau$, given on objects by
\[
H^M(V)
=
\left( M \oplus H(V), f^M_{V}=(M\xtwoheadleftarrow[]{\pi_M}  M\oplus V \xrightarrowtail[]{\id_M \oplus i_V} M\oplus H(V) ) \right)
\]
and morphisms in the obvious way, is an equivalence.
\end{Lem}

\begin{proof}
Essential surjectivity of $H^M$ follows from Proposition \ref{prop:isoReductionSplitting}. Faithfulness of $H^M$ follows from that of the hyperbolic functor $H:\mathcal{S}\rightarrow\mathcal{S}_H$. To see that $H^M$ is full, let $\phi: H^M(U) \rightarrow H^M(V)$ be a morphism, that is, an isometry $\phi:M\oplus H(U) \rightarrow M \oplus H(V) $ such that $\phi \circ f^M_U=f^M_V$ as morphisms in $ Q_{H}$. Since
\[
\phi \circ f^M_U
=
(M\xtwoheadleftarrow[]{\pi_M}  M\oplus U \xrightarrowtail[]{\phi \circ i_{M \oplus U}} M\oplus H(V)),
\]
the equality $\phi \circ f^M_U=f^M_V$ implies the existence of an isomorphism $\alpha: M \oplus U \rightarrow M \oplus V$ such that $\pi_{M \oplus V \rightarrow M} \circ \alpha = \pi_{M \oplus U \rightarrow M}$ and $(\id_M \oplus i_V) \circ \alpha=\phi \circ i_{M \oplus U}$. By the first equality, there exists a commutative diagram
\[
\begin{tikzpicture}[baseline= (a).base]
\node[scale=1] (a) at (0,0){
\begin{tikzcd}[column sep=30pt, row sep=30pt]
U \arrow[tail]{r}[above]{i_U} \arrow{d}[left]{\tilde{\alpha}} & M \oplus U \arrow[two heads]{r}[above]{\pi_M} \arrow{d}[left]{\alpha}  & M \arrow[d,equal]\\
V \arrow[tail]{r}[below]{i_V} & M \oplus V \arrow[two heads]{r}[below]{\pi_M} & M.
\end{tikzcd}
};
\end{tikzpicture}
\]
Since $\alpha$ is an isomorphism, so too is $\tilde{\alpha}$. Using Lemma \ref{lem:uniqueMorphismConflation}, we conclude that $\alpha = \id_M \oplus \tilde{\alpha}$. Using the second equality, we then verify that $\phi = \id_M \oplus H(\tilde{\alpha}) = : H^M(\tilde{\alpha})$. The $\mathcal{S}$-equivariance of $H^M$ is clear.
\end{proof}

Let $K \in \mathcal{S}$. Then $0_K:=(0\twoheadleftarrow K \yrightarrowtail[]{i_K} H(K))$ is a morphism $0 \rightarrow H(K)$ in $Q_H$. The base change functor $0_K^*:H(K)\backslash \tau\rightarrow 0\backslash\tau$ is given on objects by
\[
0_K^*(N, (H(K)\stackrel{q}{\twoheadleftarrow} E \stackrel{j}{\rightarrowtail} N))= (N, (0 \xtwoheadleftarrow[]{\quad} K\times_{H(K)}  E \xrightarrowtail[]{i_K \times_M j} N)).
\]
Consider also the functor $0_{K*}: 0\backslash \tau\rightarrow H(K)\backslash \tau$ given on objects by
$$0_{K*}(N, (0\twoheadleftarrow E \stackrel{j}{\rightarrowtail} N))=(H(K) \oplus N , (H(K) \xtwoheadleftarrow[]{\pi_{H(K)}} H(K)\oplus E \xrightarrowtail[]{\id_{H(K)} \oplus j} H(K) \oplus N))$$
and on morphisms in the obvious way.

\begin{Lem}\label{lem:homotopyInverseOfPullback}
There are natural isomorphisms of functors
\[
K\cdot\simeq0^{*}_{K}0_{K*} : 0\backslash \tau \rightarrow 0\backslash \tau
\]
and
\[
K\cdot\simeq0_{K*}0^{*}_{K}: H(K)\backslash \tau \rightarrow H(K)\backslash \tau,
\]
where $K \cdot $ denotes the action of $K \in \mathcal{S}$ on $Q_H$.
\end{Lem}

\begin{proof}
For ease of notation, set $M=H(K)$. Fix $A=(N, (0\twoheadleftarrow E \stackrel{j}{\rightarrowtail} N)) \in 0\backslash\tau$. We compute
\begin{align*} 
0^{*}_{K}0_{K*}(A)
&=\Big(M \oplus N, (0 \twoheadleftarrow K\times_M (M\oplus E) \xxrightarrowtail[]{\id_K \times_M (\id_M \oplus j)} M\oplus N)\Big)\\
&=\Big(M \oplus N, (0 \twoheadleftarrow K\oplus E \xrightarrowtail[]{i_K \oplus j} M\oplus N)\Big)\\
&= K\cdot \Big(N, (0 \twoheadleftarrow  E \yrightarrowtail[]{j} N)\Big)=K\cdot A.
\end{align*}
Lemma \ref{lem:pullbackDirectSum} was used to deduce the second equality. These identifications are clearly natural in $A$ and the first claimed natural isomorphism follows.

In view of Lemma \ref{lem:HMEquiv}, to establish the second natural isomorphism it suffices to show that $K\cdot H^M\simeq0_{K*}0^{*}_{K}H^M$. Let $V\in \mathcal{S}$ and set $N=H(V)$.
Then we have
\begin{align*} 
0_{K*}0^{*}_{K}H^M(V)
&=0_{K*}\Big(M\oplus N, (0 \twoheadleftarrow K\oplus V \xrightarrowtail[]{i_K \oplus i_V} M \oplus N)\Big)\\
&=\Big(M \oplus M\oplus N, (0 \twoheadleftarrow M\oplus K\oplus E \xxrightarrowtail[]{\id_M \oplus i_K \oplus i_V} M\oplus M\oplus N)\Big)
\intertext{while}
K\cdot H^M(V)&=\Big(M \oplus M\oplus N, (0 \twoheadleftarrow K\oplus M\oplus E \xxrightarrowtail[]{i_K \oplus \id_M \oplus i_V} M\oplus M\oplus N)\Big).
\end{align*}
The map $s_V:M \oplus M\oplus N\rightarrow M \oplus M\oplus N$, which swaps the first two summands and is the identity on the third, provides an isomorphism $s_V:K\cdot H^M(V)\rightarrow 0_{K*}0^{*}_{K}H^M(V)$ which is evidently natural in $V$.
\end{proof}

\begin{Rem}
The mistake in the argument of \cite{charney1986} is related to the swap natural transformation $s_V$. In \cite{charney1986}  it is claimed that the pullback $\sigma^*: (M\oplus M) \backslash \tau\rightarrow (M\oplus M) \backslash \tau$ associated to the swap $\sigma: M\oplus M\rightarrow M\oplus M$ is an inner automorphism (which is not true in general) and hence induces the identity on homology; see \cite[Remark 1.5]{hornbostel2004}. Our argument, which relies in a critical manner on the uniquely split assumption, uses the explicit description of the categories $M\backslash \tau$ in Lemma \ref{lem:HMEquiv} and does not require that $\sigma^*$ be inner.
\end{Rem}

The groupoid $\mathcal{S}$ acts on $\mathcal{S}_H$ and $ Q_H$ via the hyperbolic functor $H: \mathcal{S}\rightarrow \mathcal{S}_H$. With these actions, $\tau : \mathcal{S}_H \rightarrow Q_H$ is $\mathcal{S}$-equivariant. Consider the localized functor $\mathcal{S}^{-1}\tau: \mathcal{S}^{-1} \mathcal{S}_H\rightarrow \mathcal{S}^{-1}  Q_{H}$.
The natural inclusion $ Q_{H}\rightarrow \mathcal{S}^{-1}  Q_{H}$ admits a quasi-inverse $\lambda$, given on objects and morphisms by $\lambda(V, M)= M$ and
\[
\lambda(X,X\oplus V\xrightarrow[]{\alpha} V',H(X)\oplus N\xrightarrow[]{\beta} N')
=
\beta\circ(N\xtwoheadleftarrow[]{\pi_X} X\oplus N \xrightarrowtail[]{i_{X \oplus N}} H(X)\oplus N),
\]
respectively. Write $\theta: \mathcal{S}^{-1} \mathcal{S}_H \rightarrow  Q_{H}$ for the composition $\lambda \circ \mathcal{S}^{-1}\tau$.

\begin{Lem}
\label{lem:thetaBaseChangeIso}
For any morphism $\beta: M'\rightarrow M$ in $ Q_{H}$, the base change functor $\beta^*: M\backslash \theta \rightarrow M'\backslash \theta$ induces a homotopy equivalence on classifying spaces.
\end{Lem}

\begin{proof}
There is an equivalence $\Phi_M: \mathcal{S}^{-1}(M \backslash \tau) \rightarrow M \backslash \theta$; \textit{cf}. the proof of \cite[Lemma 3.2]{charney1986}. On objects, set $\Phi_M(U,(N, f)) = ((U, N), f)$. Consider a morphism
\[
(V; \gamma, \delta): (U,(N,f)) \rightarrow (U^{\prime},(N^{\prime},f^{\prime}))
\]
in $\mathcal{S}^{-1}(M \backslash \tau)$. Here $V \in \mathcal{S}$ and $\gamma: V \oplus U \xrightarrow[]{\sim} U^{\prime}$ and $\delta: H(V) \oplus N \xrightarrow[]{\sim} N^{\prime}$ is an isometry which satisfies $\beta \circ f = f^{\prime}$. Set $\Phi_M(V;\gamma,\delta) = (V; \gamma, \delta)$. It is straightforward to verify that the following diagram commutes up to natural isomorphism:
\[
\begin{tikzpicture}[baseline= (a).base]
\node[scale=1] (a) at (0,0){
\begin{tikzcd}[column sep=4em]
\mathcal{S}^{-1} (M \backslash \tau) \arrow{d}[left]{\Phi_M} \arrow{r}[above]{\mathcal{S}^{-1} (\beta^*)} & \mathcal{S}^{-1} (M^{\prime} \backslash \tau) \arrow{d}[right]{\Phi_{M^{\prime}}} \\
M \backslash \theta \arrow{r}[below]{\beta^*} & M^{\prime} \backslash \theta .
\end{tikzcd}
};
\end{tikzpicture}
\]

It therefore suffices to prove that $\mathcal{S}^{-1} (\beta^*)$ induces an equivalence on classifying spaces. Let $K^{\prime} \rightarrowtail M^{\prime}$ be a Lagrangian and set $K := K^{\prime} \times_{M^{\prime}} E \rightarrowtail M$. We then have $0_{K} \circ \beta=0_{K'}$ and, by the 2-out-of-3 property for homotopy equivalences, it suffices to consider the case $\beta = 0_K$. Lemma \ref{lem:homotopyInverseOfPullback} implies that $\mathcal{S}^{-1}(0_{K*})$ and $\mathcal{S}^{-1}(0^*_{K})$ are quasi-inverse up to $\mathcal{S}^{-1}(K\cdot)$. However, $\mathcal{S}^{-1}(K\cdot)$ is homotopic to the identity by properties of localization. The lemma follows.
\end{proof}

\begin{Lem}\label{lem:homotopyfibrationforGW}
The sequence
\[
B(\mathcal{S}^{-1}\mathcal{S}) \xrightarrow[]{B(\mathcal{S}^{-1}H)} B(\mathcal{S}^{-1}\mathcal{S}_H) \xrightarrow[]{\quad B \theta\quad} BQ_{H}
\]
is a homotopy fibration.
\end{Lem}

\begin{proof} 
In view of Lemma $\ref{lem:thetaBaseChangeIso}$, Quillen's Theorem B applies to $\theta$, allowing us to conclude that $B\theta$ is a homotopy fibration with homotopy fibre $B(0\backslash\theta)$. Lemma \ref{lem:HMEquiv} and the equivalence $\Phi_M$ from the proof of Lemma \ref{lem:thetaBaseChangeIso} compose to an equivalence
\[
\mathcal{S}^{-1} \mathcal{S} \xrightarrow[]{\mathcal{S}^{-1} H^M} \mathcal{S}^{-1}(M \backslash \tau) \xrightarrow[]{\ \ \Phi_M\ \ } M \backslash \theta.
\]
Setting $M=0$ completes the proof.
\end{proof}

\begin{proof}[Proof of Theorem \ref{thm:GWComparison}]
We follow the strategy of proof from the exact setting \cite[Theorem 3.5]{hornbostel2004}. Consider the diagram
\begin{center}
\begin{tikzcd}
                                                              & \mathcal{S}^{-1}\mathcal{E}_0 \arrow[d] \\
\mathcal{S}^{-1}\mathcal{E} \times_Q Q_H  \arrow[d, "\tilde{g}^*"'] \arrow[r] \arrow[dr, phantom, "\scalebox{1.0}{$\lrcorner$}" , very near start, color=black] & \mathcal{S}^{-1}\mathcal{E} \arrow[d, "\tilde{g}"]                     \\
Q_H \arrow[r,"F" below]                                                & Q                                                          
\end{tikzcd}
\end{center}
whose square is cartesian and right column is as in  Lemma \ref{lem:homotopyFibrationComparisonKTheory}. Since $\tilde{g}$ is fibred and fulfills the conditions of Quillen's Theorem B and $\mathcal{S}^{-1}\mathcal{E}$ is contractible, we can apply \cite[Proposition 3.4]{hornbostel2004} to conclude that (after passing to classifying spaces) the sequences
\[
\mathcal{S}^{-1}\mathcal{E}_0 \xrightarrow[]{\iota}\mathcal{S}^{-1}\mathcal{E} \times_Q Q_H \xrightarrow[]{\tilde{g}^*} Q_H
\]
and
\[
\mathcal{S}^{-1}\mathcal{E} \times_Q Q_H \xrightarrow[]{\tilde{g}^*} Q_H \xrightarrow[]{F} Q
\]
are homotopy fibrations. Consider the diagram
\begin{center}
\begin{tikzcd}
\mathcal{S}^{-1}\mathcal{S}  \arrow[r, "\mathcal{S}^{-1}H"] \arrow[d,"\mathcal{S}^{-1}F_0",swap]& \mathcal{S}^{-1}\mathcal{S}_H \arrow[r, "\theta"] \arrow[d,"J"] & Q_H   \arrow[d,equal]                                                  \\
\mathcal{S}^{-1}\mathcal{E}_0 \arrow[r, "\iota" below] &\mathcal{S}^{-1}\mathcal{E} \times_Q Q_H \arrow[r, "\tilde{g}^*" below] & Q_H,
\end{tikzcd}
\end{center}
where $F_0$ is the equivalence from Lemma \ref{lem:descriptionoffibers} and $J$ is defined as in \cite[Theorem 3.5]{hornbostel2004}. The right square commutes and the left square commutes up to natural isomorphism. By Lemma \ref{lem:homotopyfibrationforGW}, the top row is a homotopy fibration. Since $\mathcal{S}^{-1} F_0$ is an equivalence, it follows from the five lemma that $BJ$ is a weak homotopy equivalence.

The above discussion gives a weak homotopy equivalence
\[
\mathcal{GW}_H \simeq B(\mathcal{S}^{-1}\mathcal{E}) \times_{BQ} BQ_H.
\]
Moreover, there is a homotopy equivalence $B(\mathcal{S}^{-1}\mathcal{S}_H) \simeq B(\mathcal{S}_H^{-1}\mathcal{S}_H) = \mathcal{GW}_H^\oplus$. Composing with $BJ$ then gives the desired weak homotopy equivalence $\mathcal{GW}_H \simeq \mathcal{GW}^{\oplus}_H$.
\end{proof}

 The above arguments show that there is a homotopy commutative diagram
\[
\begin{tikzpicture}[baseline= (a).base]
\node[scale=1] (a) at (0,0){
\begin{tikzcd}[column sep=4.0em]
B(\mathcal{S}^{-1} \mathcal{S}) \arrow{d} \arrow{r}[above]{B(\mathcal{S}^{-1}H)} & B(\mathcal{S}^{-1} \mathcal{S}_H) \arrow{d} \arrow{r}[above]{B \theta} & BQ_H(\mathcal{A}) \arrow[equal]{d}\arrow{r}[above]{BF} & BQ(\mathcal{A}) \arrow[equal]{d}\\
\Omega BQ(\mathcal{A}) \arrow[equal]{d} \arrow{r} & \mathcal{GW}_H(\mathcal{A}) \arrow{d} \arrow{r} & BQ_H(\mathcal{A}) \arrow{d} \arrow{r}[above]{BF} & BQ(\mathcal{A}) \arrow[equal]{d}\\
\Omega BQ(\mathcal{A}) \arrow{r}[below]{H_*}& \mathcal{GW}(\mathcal{A}) \arrow{r} & BQ_h(\mathcal{A}) \arrow{r}[below]{BF} & BQ(\mathcal{A})
\end{tikzcd}
};
\end{tikzpicture}
\]
whose rows are homotopy fibre sequences. Using Theorem \ref{thm:KComparison} we conclude that the first two rows are weak homotopy equivalent. It follows that the morphism $H_{*}$ from equation \eqref{eq:wittCoker} is explicitly realized by the hyperbolic functor $H$. Moreover, since $BQ_H(\mathcal{A}) \subset BQ_h(\mathcal{A})$ is the connected component of $0$, we obtain from Theorem \ref{thm:GWComparison} the following result

\begin{Cor}
\label{cor:higherGWGroupsQ+}
For each $i \geq 1$, there is an isomorphism of abelian groups
\[
GW_i(\mathcal{A}) 
\simeq
GW^{\oplus}_{H,i}(\mathcal{A}).
\]
\end{Cor}

It remains to compute $GW_0(\mathcal{A})$ and to determine the homotopy type of the connected components of $\mathcal{GW}(\mathcal{A})$ which do not contain $0$. This is done on the next section. 
\subsection{The connected components of $\mathcal{GW}(\mathcal{A})$}
\label{sec:connCptGW}

We continue to denote by $\mathcal{A}$ a uniquely split proto-exact category with duality. In the previous section, we determined the homotopy groups $GW_i(\mathcal{A})$, $i \geq 1$, in terms of $\mathcal{GW}_H^{\oplus}(\mathcal{A})$. However, we did not describe $GW_0(\mathcal{A})$. As discussed in Section \ref{sec:GWProtoEx}, both $GW_0(\mathcal{A})$ and $W_0(\mathcal{A})$ are (commutative) monoids, and not groups, due to the lack of a diagonal map for $\oplus$. This complicates their computation. To get around this problem, we henceforth restrict attention to combinatorial categories.

Denote by $\mathcal{S}_{I}\subset \mathcal{S}_h$ the full subcategory of isotropically simple symmetric forms. Given $S\in\mathcal{S}_{I}$, we abbreviate $\operatorname{Aut}_{\mathcal{S}_h}(S)$ to $G_S$ and write $\mathbf{B} G_S$ for the associated groupoid with one object. The geometric realization of $\mathbf{B} G_S$ is $BG_S$.

\begin{Lem}
\label{lem:isoSimpProperties}
Let $\mathcal{A}$ be a uniquely split combinatorial noetherian proto-exact category with duality.
\begin{enumerate}[wide,labelwidth=!, labelindent=0pt,label=(\roman*)]
\item \label{lemPart:stableDirSum}The category $\mathcal{S}_{I}$ is closed under direct sum. In particular, the set $\pi_0(\mathcal{S}_I)$ of connected components of $\mathcal{S}_I$ is a commutative monoid.

\item \label{lemPart:WittIso} The natural map $\pi_0(\mathcal{S}_I) \rightarrow W_0(\mathcal{A})$ is an isomorphism of monoids.

\item Direct sum induces an equivalence of categories $\mathcal{S}_I \times \mathcal{S}_H \rightarrow \mathcal{S}_h$.
\end{enumerate}
\end{Lem}

\begin{proof}
Let $S_1, S_2 \in \mathcal{S}_I$ and suppose that $U \rightarrowtail S_1 \oplus S_2$ is isotropic. The combinatorial assumption implies that $U \cap S_i \rightarrowtail S_i$, $i=1,2$, is isotropic and hence each $U \cap S_i$ is the zero object. It follows that $U \simeq 0$ and $S_1 \oplus S_2$ is isotropically simple.

By definition, $W_0(\mathcal{A})$ is the quotient of $\pi_0(\mathcal{S}_h)$ by the submonoid of metabolic objects. By Lemma \ref{lem:noMeta}, metabolics are necessarily hyperbolic. The second statement therefore follows from Proposition \ref{prop:isoSimpDecomp}.

The third statement also follows from Proposition \ref{prop:isoSimpDecomp}: essentially surjectivity follows from part \ref{propPart:essSurj} and fully faithfulness follows from the proof of part \ref{propPart:fullFaith}.
\end{proof}

The remainder of this section is devoted to proving the following theorem.

\begin{Thm}
\label{thm:GWComputation}
Let $\mathcal{A}$ be a uniquely split combinatorial noetherian proto-exact category with duality. Then there is a weak homotopy equivalence
\[
\mathcal{GW}(\mathcal{A})
\simeq
\bigsqcup_{w \in W_0(\mathcal{A})}BG_{S_w}\times\mathcal{GW}_H(\mathcal{A}),
\]
where $S_w \in \mathcal{S}_I$ is an isotropically simple representative of $w \in W_0(\mathcal{A})$ under the isomorphism of Lemma \ref{lem:isoSimpProperties}\ref{lemPart:WittIso}.
\end{Thm}

For each $S\in \mathcal{S}_I$, denote by $Q^S \subset Q_h$ the full subcategory of objects isometric to an object of the form $S\oplus H(V)$ for some $V \in \mathcal{A}$. Let $\mathcal{GW}^S$ be the homotopy fibre of $BF: BQ^S\rightarrow BQ$ over $0$. With this notation, we have $Q^0=Q_H$ and $\mathcal{GW}^0=\mathcal{GW}_H$.

\begin{Lem}
\label{lem:QandGWDecomp}
Work in the setting of Theorem \ref{thm:GWComputation}.
\begin{enumerate}[wide,labelwidth=!, labelindent=0pt,label=(\roman*)]
\item \label{eq:splittingQS} For each $S \in \mathcal{S}_I$, there is an equivalence of categories $Q^S \simeq \mathbf{B} G_S \times Q_H$.

\item The category $Q_h$ decomposes into connected components as
$$Q_h=\bigsqcup_{w \in W_0(\mathcal{A})} Q^{S_w}.$$
\item The space $\mathcal{GW}$ decomposes into (not-necessarily connected) components as
\[
\mathcal{GW}=\bigsqcup_{w \in W_0(\mathcal{A})} \mathcal{GW}^{S_w}.
\]
\end{enumerate}
\end{Lem}

\begin{proof}
Let $\eta: BG_S \times Q_H \rightarrow Q^S$ be the functor given on objects and morphisms by $H(V) \mapsto S \oplus H(V)$ and
\[
\left(\phi, H(V) \xtwoheadleftarrow[]{\!\!\!\pi\!\!\!} E \yrightarrowtail[]{i} H(V^{\prime}) \right)
\mapsto
\left(
S \oplus H(V) \xtwoheadleftarrow[]{\!\!\!\id_S \oplus \pi\!\!\!} S \oplus E \xrightarrowtail[]{\!\!\!\phi \oplus i\!\!\!} S \oplus H(V^{\prime})
\right),
\]
respectively. We claim that $\eta$ is an equivalence. Essentially surjectivity of $\eta$ follows from the definition of $Q^S$. Because $S$ is isotropically simple, the combinatorial assumption ensures that any isotropic morphism with target $S \oplus H(V^{\prime})$ factors through $i_{H(V^{\prime})} : H(V^{\prime}) \rightarrowtail S \oplus H(V^{\prime})$. It follows from this that $\eta$ is full. A short matrix calculation shows that $\eta$ is faithful.

Consider the second statement. By Proposition \ref{prop:isoSimpDecomp}, the canonical functor $\bigsqcup_{w \in W_0} Q^{S_w} \rightarrow Q_h$ is essentially surjective. Suppose that there is a morphism $M_1 \rightarrow M_2$ in $Q_h$. By Propositions \ref{prop:isoReductionSplitting} and \ref{prop:isoSimpDecomp}, there exist isometries $M_1 \simeq S \oplus H(U)$ and $M_2 \simeq S \oplus H(V)$ for a unique (up to isometry) $S \in \mathcal{S}_I$. It follows from Lemma \ref{lem:isoSimpProperties}\ref{lemPart:WittIso} that $M_1$ and $M_2$ represent the same class in $W_0(\mathcal{A})$, proving fully faithfulness of $\bigsqcup_{w \in W_0} Q^{S_w} \rightarrow Q_h$. Since $S_{w} \in Q^{S_w}$, the category $Q^{S_w}$ is connected.

Since $\mathcal{GW}$ is a homotopy fibre of $BQ_h \rightarrow BQ$, the third statement follows from the second. 
\end{proof}

By Lemma \ref{lem:QandGWDecomp}, in order to understand $\mathcal{GW}$, it suffices to describe the components $\mathcal{GW}^S$. 
Similarly to the argument in the proof of Theorem \ref{thm:GWComparison}, which treats the case $S=0$, we can apply \cite[Proposition 3.4]{hornbostel2004} to the square
\begin{center}
\begin{tikzcd}
\mathcal{S}^{-1}\mathcal{E} \times_Q Q^{S}  \arrow[d, "\tilde{g}^*"'] \arrow[r] \arrow[dr, phantom, "\scalebox{1.0}{$\lrcorner$}" , very near start, color=black] & \mathcal{S}^{-1}\mathcal{E} \arrow[d, "\tilde{g}"]                     \\
Q^{S} \arrow[r,"F" below]                                                & Q                                                           
\end{tikzcd}
\end{center}
to see that $\mathcal{GW}^S$ is homotopy equivalent to the geometric realization of the category
$$\mathcal{S}^{-1}\mathcal{E} \times_Q Q^S.$$
We therefore proceed to consider the category $\mathcal{S}^{-1}\mathcal{E} \times_Q Q^S.$

\begin{Prop}
\label{prop:DecompositionOfModelOfConnectectComponentOfGW}
For each $S \in \mathcal{S}_I$, there is an equivalence of categories
$$BG_S\times (\mathcal{S}^{-1}\mathcal{E} \times_Q Q_H)\simeq \mathcal{S}^{-1}\mathcal{E} \times_Q Q^S.$$
\end{Prop}

\begin{proof}
Consider the category $\mathcal{E} \times_Q Q^S$. By definition, an object is a triple
$$((A\rightarrowtail B \twoheadrightarrow C)\in \mathcal{E}, M\in Q^S, (\phi: C \xrightarrow[]{\sim} M) \in Q)$$
and a morphism is a pair of morphisms in $\mathcal{E}$ and $Q^S$ which make the obvious square in $Q$ commute. We claim that there is an equivalence
\[
\mathcal{E} \times_Q Q^S \simeq BG_{S}\times \mathcal{E} \times_Q Q_{H}.
\] 
Using Propositions \ref{prop:isoReductionSplitting} and \ref{prop:isoSimpDecomp}, we see that $\mathcal{E} \times_Q Q^S$ is equivalent to the category 
 $\mathcal{E}^{S}$ whose objects are conflations of the form
\begin{center}
\begin{tikzcd}
A \arrow["i_{A}",r, tail] & S\oplus H\oplus A  \arrow["\pi_{S\oplus H}",r, two heads] &  S\oplus H,         
\end{tikzcd}
\end{center}
for some hyperbolic object $H$, and whose morphisms are equivalence classes of commutative diagrams of the form 
\begin{center}
\begin{tikzcd}[column sep=4em]
A \arrow["i_{A}",r, tail]                                                 & S\oplus H \oplus A \arrow["\pi_{S\oplus H}",r, two heads] \arrow[d, equal] &  S\oplus H                               \\
A' 
\arrow["\alpha",u, tail] 
\arrow[d, equal] 
\arrow["i_{A} \circ \alpha",r, tail] 
& S\oplus H \oplus A\arrow["\beta",d, tail] \arrow["\pi",r, two heads]               & 
S\oplus W 
\arrow["\id_{S} \oplus p",swap,u, two heads] 
\arrow["\phi \oplus i"',swap, d, tail] \\
A' \arrow["i_{A'}",swap, r, tail]                                                 & S\oplus H' \oplus A' \arrow["\pi_{S'\oplus H'}",swap, r, two heads]                                & S\oplus H'    .                                 
\end{tikzcd}
\end{center}
In the previous diagram, all rows are conflations and the right column is a morphism in $Q^S$ or, by Lemma \ref{lem:QandGWDecomp}\ref{eq:splittingQS}, in $BG_S \times Q_H$. By splitting the morphism $p,$ any such diagram is equivalent to one of the form
\begin{equation}
\label{diag:bigSquare}
\begin{tikzcd}[column sep=4em]
A \arrow["i_{A}",r, tail]                                                 & S\oplus H \oplus A \arrow["\pi_{S\oplus H}",r, two heads] \arrow[d, equal] &  S\oplus H                               \\
A' 
\arrow["\alpha",u, tail] 
\arrow[d, equal] 
\arrow["i_{A} \circ \alpha",r, tail] 
& S\oplus H \oplus A\arrow["\beta",d, tail] \arrow["\pi",r, two heads]               & 
S\oplus H\oplus X 
\arrow["\pi_{S\oplus H}", swap,u, two heads] 
\arrow["\phi \oplus i"', swap,d, tail] \\
A' \arrow["i_{A'}", swap,r, tail]                                                 & S\oplus H' \oplus A' \arrow["\pi_{S \oplus H'}", swap,r, two heads]                                & S\oplus H'                                      
\end{tikzcd}
\end{equation}
where $W\simeq H\oplus X$.  By the commutativity of the top right square, $\pi$ is of the form
\[
\left( \begin{array}{rrrr}
\id_{S} & 0 & 0 \\
0 & * & *  \\
x & * & * \\
\end{array}\right): S\oplus H\oplus A \twoheadrightarrow S\oplus H\oplus X.
\]
By the commutativity of the bottom right square, $\beta$ is of the form
\[
\beta=\left( \begin{array}{rrrr}
\phi & 0 & 0 \\
0 & * & *  \\
z & * & * \\
\end{array}\right): S\oplus H\oplus A \rightarrowtail S\oplus H'\oplus A'.
\]
Since $\phi$ is an isomorphism, Lemma \ref{lem:partialIsoComb} implies that $z=0$. Commutativity of right hand squares then gives $x=\pi_X \circ \pi \circ i_S=0$ and $\pi_{H\oplus X} \circ \pi \circ i_S=0$, with the latter equality following from the equalities
\[
(\phi\oplus i) \circ i_{H\oplus X} \circ \pi_{H\oplus X} \circ \pi \circ i_S=i_{H'} \circ \pi_{H'} \circ (\phi \oplus i) \circ \pi \circ i_S=\pi_{S\oplus H'} \circ \beta \circ i_S=0
\]
and the fact that $(\phi\oplus i) \circ i_{H\oplus X}$ is a monomorphism. Hence, $S$ decouples from the diagram \eqref{diag:bigSquare} and the obvious functor $BG_{S}\times \mathcal{E}^{0} \to \mathcal{E}^{S}$ is an equivalence. By the argument above, specialized to the case $S=0$, there is an equivalence $ \mathcal{E}^{0}\simeq \mathcal{E}\times_{Q}Q_{H}$. Summarizing, we obtain an equivalence $BG_{S}\times \mathcal{E}\times_{Q}Q_{H}\simeq \mathcal{E}\times_{Q}Q^{S}$.

Let $\mathcal{S}$ act on $\mathcal{E} \times_Q Q^S$ by the standard action on $\mathcal{E}$ (see Section \ref{sec:Q=+Thm}) and on $Q^S$ by the trivial action. Then there is an equivalence $\mathcal{S}^{-1}\mathcal{E} \times_Q Q^S \simeq \mathcal{S}^{-1}(\mathcal{E} \times_Q Q^S)$ given by simply reordering the bracketing of objects and morphisms. Hence, there are  equivalences 
\[
BG_{S}\times (\mathcal{S}^{-1}\mathcal{E}\times_{Q}Q_{H})\simeq BG_{S}\times \mathcal{S}^{-1}(\mathcal{E}\times_{Q}Q_{H})\simeq\mathcal{S}^{-1}(\mathcal{E}\times_{Q}Q^{S})\simeq\mathcal{S}^{-1}\mathcal{E}\times_{Q}Q^{S}.
\]
This completes the proof.
\end{proof}

\begin{proof}[Proof of Theorem \ref{thm:GWComputation}]
The space $\mathcal{GW}^S$ is homotopy equivalent to the geometric realization of $\mathcal{S}^{-1}\mathcal{E}\times_QQ^S$. By Proposition \ref{prop:DecompositionOfModelOfConnectectComponentOfGW}, there is an equivalence
\[
\mathcal{S}^{-1}\mathcal{E} \times_Q Q^S \simeq BG_S\times (\mathcal{S}^{-1}\mathcal{E} \times_Q Q_H).
\]
Since the geometric realization of $\mathcal{S}^{-1}\mathcal{E} \times_Q Q_H$ is weak homotopy equivalent to $\mathcal{GW}_H$, this completes the proof.
\end{proof}

\begin{Cor}
\label{cor:GW0Description}
In the setting of Theorem \ref{thm:GWComputation}, there is an isomorphism of monoids
\[
GW_0(\mathcal{A})
\simeq
W_0(\mathcal{A}) \times GW_{H,0}(\mathcal{A}).
\]
\end{Cor}
\begin{proof}
By Proposition \ref{prop:DecompositionOfModelOfConnectectComponentOfGW}, we see that $\mathcal{GW}^S(\mathcal{A}) \simeq BG_S \times \mathcal{GW}_H(\mathcal{A})$ for each $S \in \mathcal{S}_I$. In particular, there is a monoid isomorphism $\pi_0 \mathcal{GW}^S(\mathcal{A}) \simeq GW_{H,0}(\mathcal{A})$. In view of Theorem \ref{thm:GWComputation}, it follows that the stated isomorphism holds at the level of sets. Using Lemma \ref{lem:isoSimpProperties}\ref{lemPart:stableDirSum}, it is easy to see that the decomposition in Theorem \ref{thm:GWComputation} is compatible with direct sum, which implies the statement for monoids.
\end{proof}

Using Corollary \ref{cor:GW0Description} and the discussion above Corollary \ref{cor:higherGWGroupsQ+} to identify $\pi_0 H_*$ with the hyperbolic map, we find a monoid isomorphism
\[
W_0(\mathcal{A})
\simeq
\coker(K_0(\mathcal{A}) \xrightarrow[]{\pi_0 H_*} GW_0(\mathcal{A})).
\]
In particular, the definition of the higher Witt groups $W_i(\mathcal{A})$, $i \geq 1$, given in Section \ref{sec:hermQConstr} extends compatibly to $i=0$.

\begin{ex}\label{ex: GW-theory of Vect F1}
We continue our running test case $\Vect_\Fun$ from Examples \ref{ex: Vect F1 as proto-exact category}, \ref{ex: duality for Vect F1} and \ref{ex: K-theory of Vect F1}. As shown in \cite[\S 2.3]{eberhardtLorscheidYoung2020b}, we have homotopy equivalences
\[
\mathcal{GW}_H(\Vect_\Fun) \simeq \mathcal{GW}_H^\oplus(\Vect_\Fun)
\simeq
\mathbb{Z} \times B \big((\mathbb{Z} \slash 2) \wr \Sigma_{\infty}\big)^+
\]
and
\[
\mathcal{GW}^{\oplus}(\Vect_{\mathbb{F}_1}) \simeq \mathbb{Z}^2 \times B\left( \big( (\mathbb{Z} \slash 2) \wr \Sigma_{\infty} \big) \times \Sigma_{\infty} \right)^+
\]
and a weak homotopy equivalence
\[
\mathcal{GW}(\Vect_{\mathbb{F}_1}) \simeq \bigsqcup_{n \in \mathbb{Z}_{\geq 0}}B \Sigma_n \times \mathbb{Z} \times B ((\mathbb{Z} \slash 2) \wr \Sigma_{\infty})^+.
\]
This exemplifies, in particular, that the $Q$-construction and the $+$-construction lead to significantly different Grothendieck--Witt spaces for proto-exact categories.
\end{ex}

\bibliographystyle{plain}
\bibliography{mybib}

\def\cprime{$'$} \def\cprime{$'$}
\begin{thebibliography}{10}

\bibitem{charney1986}
R.~Charney and R.~Lee.
\newblock On a theorem of {G}iffen.
\newblock {\em Michigan Math. J.}, 33(2):169--186, 1986.

\bibitem{chu2012}
C.~Chu, O.~Lorscheid, and R.~Santhanam.
\newblock Sheaves and {$K$}-theory for {$\Bbb F_1$}-schemes.
\newblock {\em Adv. Math.}, 229(4):2239--2286, 2012.

\bibitem{chu2010}
C.~Chu and J.~Morava.
\newblock On the algebraic {$K$}-theory of monoids.
\newblock ar{X}iv:1009.3235, 2010.

\bibitem{connes2010}
A.~Connes and C.~Consani.
\newblock Schemes over {$\Bbb F_1$} and zeta functions.
\newblock {\em Compos. Math.}, 146(6):1383--1415, 2010.

\bibitem{deitmar2006}
A.~Deitmar.
\newblock Remarks on zeta functions and {$K$}-theory over {${\bf F}_1$}.
\newblock {\em Proc. Japan Acad. Ser. A Math. Sci.}, 82(8):141--146, 2006.

\bibitem{deitmar2012}
A.~Deitmar.
\newblock Belian categories.
\newblock {\em Far East J. Math. Sci. (FJMS)}, 70(1):1--46, 2012.

\bibitem{dyckerhoff2019}
T.~Dyckerhoff and M.~Kapranov.
\newblock {\em Higher {S}egal Spaces}.
\newblock Lecture Notes in Mathematics. Springer, 2019.

\bibitem{eberhardtLorscheidYoung2020b}
J.~Eberhardt, O.~Lorscheid, and M.~Young.
\newblock Algebraic {$K$}-theory and {G}rothendieck--{W}itt theory of monoid
  schemes.
\newblock ar{X}iv:2009.XXXXX, 2020.

\bibitem{eppolito2020}
C.~Eppolito, J.~Jun, and M.~Szczesny.
\newblock Proto-exact categories of matroids, {H}all algebras, and
  {$K$}-theory.
\newblock {\em Math. Z.}, 296(1-2):147--167, 2020.

\bibitem{grayson1976}
D.~Grayson.
\newblock Higher algebraic {$K$}-theory. {II} (after {D}aniel {Q}uillen).
\newblock In {\em Algebraic {$K$}-theory ({P}roc. {C}onf., {N}orthwestern
  {U}niv., {E}vanston, {I}ll., 1976)}, pages 217--240. Lecture Notes in Math.,
  Vol. 551, 1976.

\bibitem{hekking2017}
J.~Hekking.
\newblock {\em Segal objects in homotopical categories \& {$K$}-theory of
  proto-exact categories}.
\newblock 2017.
\newblock Thesis (Master's)--Mathematisch Instituut, Universiteit Leiden.

\bibitem{hornbostel2002}
J.~Hornbostel.
\newblock Constructions and d\'{e}vissage in {H}ermitian {$K$}-theory.
\newblock {\em $K$-Theory}, 26(2):139--170, 2002.

\bibitem{hornbostel2004}
J.~Hornbostel and M.~Schlichting.
\newblock Localization in {H}ermitian {$K$}-theory of rings.
\newblock {\em J. London Math. Soc. (2)}, 70(1):77--124, 2004.

\bibitem{quillen1973}
D.~Quillen.
\newblock Higher algebraic {$K$}-theory. {I}.
\newblock In {\em Algebraic {$K$}-theory, {I}: {H}igher {$K$}-theories ({P}roc.
  {C}onf., {B}attelle {M}emorial {I}nst., {S}eattle, {W}ash., 1972)}, pages
  85--147. Lecture Notes in Math., Vol. 341. Springer, Berlin, 1973.

\bibitem{schlichting2004}
M.~Schlichting.
\newblock Hermitian {$K$}-theory on a theorem of {G}iffen.
\newblock {\em $K$-Theory}, 32(3):253--267, 2004.

\bibitem{schlichting2010}
M.~Schlichting.
\newblock Hermitian {$K$}-theory of exact categories.
\newblock {\em J. K-Theory}, 5(1):105--165, 2010.

\bibitem{schlichting2017}
M.~Schlichting.
\newblock Hermitian {$K$}-theory, derived equivalences and {K}aroubi's
  fundamental theorem.
\newblock {\em J. Pure Appl. Algebra}, 221(7):1729--1844, 2017.

\bibitem{schlichting2019}
M.~Schlichting.
\newblock Higher {$K$}-theory of forms {I}: From rings to exact categories.
\newblock {\em J. Inst. Math. Jussieu}, pages 1--69, 2019.

\bibitem{srinivas2008}
V.~Srinivas.
\newblock {\em Algebraic {$K$}-theory}.
\newblock Modern Birkh\"{a}user Classics. Birkh\"{a}user Boston, Inc., Boston,
  MA, second edition, 2008.

\bibitem{uridia1990}
M.~Uridia.
\newblock {$U$}-theory of exact categories.
\newblock In {\em {$K$}-theory and homological algebra ({T}bilisi, 1987--88)},
  volume 1437 of {\em Lecture Notes in Math.}, pages 303--313. Springer,
  Berlin, 1990.

\bibitem{waldhausen1985}
F.~Waldhausen.
\newblock Algebraic {$K$}-theory of spaces.
\newblock In {\em Algebraic and geometric topology ({N}ew {B}runswick,
  {N}.{J}., 1983)}, volume 1126 of {\em Lecture Notes in Math.}, pages
  318--419. Springer, Berlin, 1985.

\bibitem{weibel2013}
C.~Weibel.
\newblock {\em The {$K$}-book: An introduction to algebraic $K$-theory}, volume
  145 of {\em Graduate Studies in Mathematics}.
\newblock American Mathematical Society, Providence, RI, 2013.

\bibitem{mbyoung2018b}
M.~Young.
\newblock Relative 2-{S}egal spaces.
\newblock {\em Algebr. Geom. Topol.}, 18(2):975--1039, 2018.

\bibitem{mbyoung2021}
M.~Young.
\newblock Degenerate versions of {G}reen's theorem for {H}all modules.
\newblock {\em J. Pure Appl. Algebra}, 225(4):106557, 2021.

\end{thebibliography}
 

\end{document}